\documentclass[graybox]{svmult}

\usepackage{mathptmx}       
\usepackage{helvet}         
\usepackage{courier}        
\usepackage{type1cm}        
\usepackage{makeidx}         
\usepackage{graphicx}        
\usepackage{multicol}        
\usepackage[bottom]{footmisc}
\usepackage{amsmath,amssymb,amsfonts}
\usepackage{caption}

\makeindex             

\newcommand{\N}{{\mathbb{N}}}

\newcommand{\R}{{\mathbb{R}}}
\newcommand{\C}{{\mathbb{C}}}

\begin{document}

\title*{Optimization-based Motion Planning in Virtual Driving Scenarios with Application to Communicating Autonomous Vehicles }
\titlerunning{Optimization-based Motion Planning}
\author{Matthias Gerdts and Bj\"orn Martens}
\institute{Matthias Gerdts \at Institute of Mathematics and Applied Computing, 
  Department of Aerospace Engineering, University of the Federal Armed Forces at Munich, 
  Werner-Heisenberg-Weg 39, 85577 Neubiberg, Germany, \email{matthias.gerdts@unibw.de}
\and Bj\"orn Martens \at Institute of Mathematics and Applied Computing, 
  Department of Aerospace Engineering, University of the Federal Armed Forces at Munich, 
  Werner-Heisenberg-Weg 39, 85577 Neubiberg, Germany, \email{bjoern.martens@unibw.de}}
%
%
\maketitle

\abstract{The paper addresses the problem of providing suitable reference trajectories 
in motion planning problems for autonomous vehicles. Among the various approaches to compute 
a reference trajectory, our aim is to find those trajectories which optimize a given performance 
criterion, for instance fuel consumption, comfort, safety, time, and obey constraints, e.g.  
collision avoidance, safety regions, control bounds. This task can be approached by geometric 
shortest path problems or by optimal control problems, which need to be solved efficiently. 
To this end we use direct discretization schemes and model-predictive control in combination with 
sensitivity updates to predict optimal solutions in the presence of perturbations. 
Applications arising in autonomous driving are presented. In particular, a distributed 
control algorithm for traffic scenarios with several autonomous vehicles that use car-to-car 
communication is introduced. }

\section{Introduction}
\label{sec:intro}

Virtual driving summarizes the ability to simulate vehicles in different driving scenarios on a computer. 
It is an important tool as it allows to analyze the dynamic behavior of a vehicle and the performance of 
driver assistance systems in parallel to the development process. 
The employment of virtual driving allows to reduce costs since simulations are less expensive and less 
time consuming than real test-drives. Moreover, virtual driving is particularly useful in scenarios that 
are potentially dangerous for human drivers such as collision avoidance scenarios or driving close to 
the physical limit. Moreover, the autonomous driving strategies can be developed and analyzed in virtual 
driving systems. However, virtual driving methods cannot fully replace physical testing since 
the virtual system is based on modeling assumptions that need to be verified in practice. 
Virtual driving simulators require models for the vehicle, the road, the environment (i.e. other 
cars, pedestrians, obstacles, ...), components (i.e. sensors, cameras, ...), and the driver 
(i.e. path planning, controller, driver assistance, ...). 

In this paper we focus on the driver, suitable path planning strategies, and control actions. 
Automatic path planning strategies are in the core of every virtual or real driving scenario 
for autonomous vehicles. Amongst the various approaches, e.g. sampling methods,  
shortest path problems and optimal control techniques, we focus on deterministic 
optimization-based methods such as geometric shortest paths, optimal control, and 
model-predictive control. While in virtual driving real-time capability is only of minor interest, 
in online computations it is the most important issue. In both cases robust methods are required 
that are capable of providing a suitable result reliably.  
A discussion of technical, legal, and social aspects of autonomous driving can be found in 
the recent book \cite{Maurer2015}.

This paper is organized as follows: Section~\ref{sec:modeling} introduces working models 
for the vehicle, the road, obstacles, and the driver. In Section~\ref{sec:trajectory} approaches 
for optimization-based path planning are discussed. Amongst them we discuss geometric shortest 
path problems and optimal control approaches in more detail. 
Section~\ref{Sec:CollisionDetection} addresses collision detection methods. 
A model-predictive control scheme for communicating vehicles is suggested in 
Section~\ref{sec:control}. In Section~\ref{sec:Tracking} we discuss a feedback controller 
based on inverse kinematics for tracking a reference spline curve.

\section{Modeling} 
\label{sec:modeling}

Virtual driving requires a sufficiently realistic vehicle model, 
a model for its environment, and a driver model. Vehicle models exist in various levels of 
accuracy, ranging from simple point-mass models through single-track models to full car models. 
Which model to use depends on the effects that one likes to investigate. For handling purposes 
and online computations often a simple point-mass model or a single-track model with realistic 
tyre characteristics, compare \cite{Ger05a}, are sufficiently accurate. For the investigation of 
dynamic load changes or vibrations a full car model in terms of a mechanical multi-body system 
is necessary, compare \cite{Ger03b,Rill2011} for full car models and \cite{Burger2010,Burger2013} for 
applications with load excitations. Often, very detailed component models such as tyre models 
become necessary to investigate tyre-road contacts, compare \cite{Gallrein2007,Gallrein2014}.

\begin{figure}
  \begin{center}
    \includegraphics[scale=0.8]{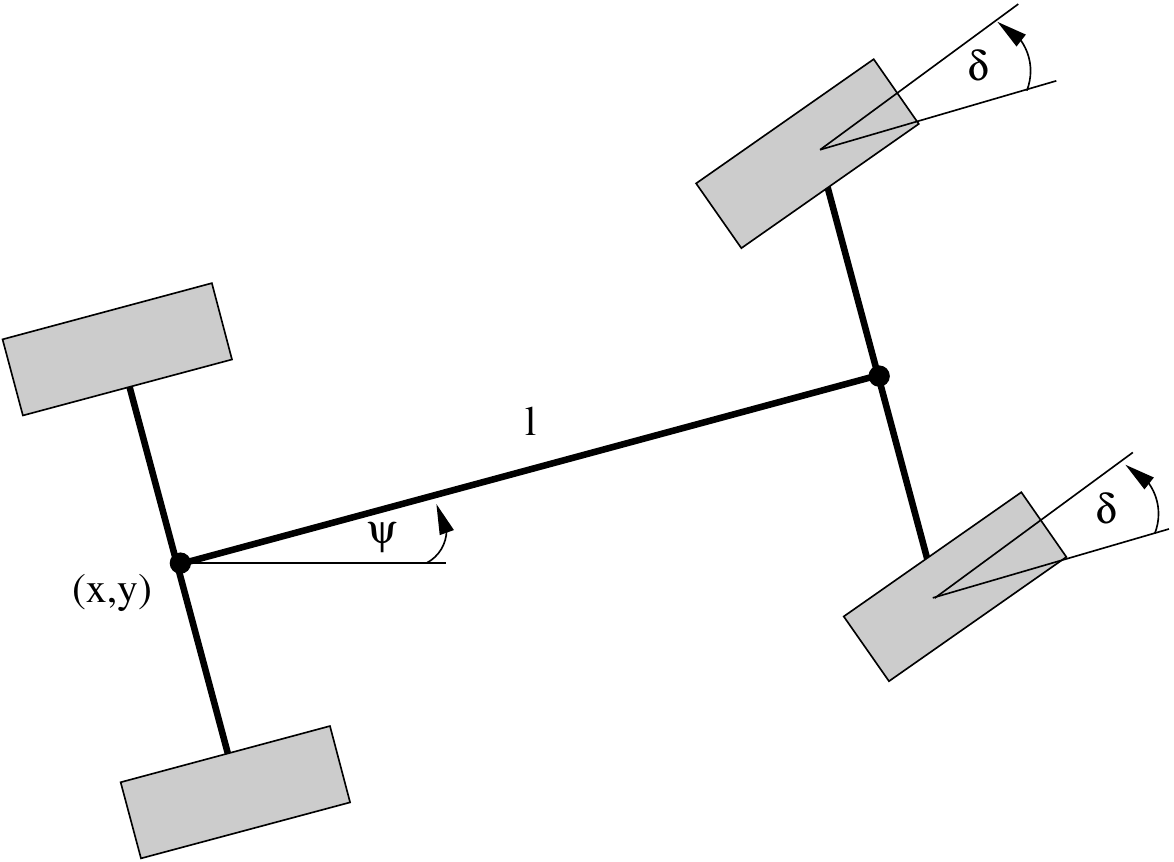}
  \end{center}
  \caption{Configuration of the kinematic car model.}\label{Fig:1}
\end{figure}

Throughout the paper we use a simple kinematic car model, compare \cite{Rill2011}, 
since it is sufficient to introduce the basic ideas. The kinematic car model is valid for low lateral 
accelerations, low lateral tyre forces, and negligible side slip angles. It is 
not suitable for investigations close to the dynamic limit, though. 
The configuration of the car in a reference coordinate system is depicted in 
Figure~\ref{Fig:1}. Herein, $\delta$ 
denotes the steering angle at the front wheels, $v$ the velocity of the car, 
$\psi$ the yaw angle, $\ell$ the distance from rear axle to front axle and 
$(x,y)$ the position of the midpoint of the rear axle. 

\medskip
The equations of motion are derived as follows. 
The midpoint position and velocity of the rear axle of the car in the fixed 
reference system compute to
\begin{eqnarray*}
  r_R = \left(\begin{array}{c} x \\  y \end{array}\right),\qquad
  v_R = \left(\begin{array}{c} x' \\  y' \end{array}\right).
\end{eqnarray*}
The midpoint position and velocity of the front axle of the car in the fixed 
reference system compute to
\begin{eqnarray*}
  r_F & = & \left(\begin{array}{c} x_F \\  y_F \end{array}\right) = r_R + S(\psi)
  \left(\begin{array}{c} \ell \\ 0 \end{array}\right) = 
  \left(\begin{array}{c} x + \ell \cos\psi \\  y + \ell
      \sin\psi \end{array}\right), \\
  v_F & = & \left(\begin{array}{c} x_F' \\  y_F' \end{array}\right) = 
  \left(\begin{array}{c} x' - \ell \psi' \sin\psi \\  y' + \ell
      \psi' \cos\psi \end{array}\right),
\end{eqnarray*}
where 
\begin{equation}\label{EQ:Drehmatrix}
  S(\psi) = \left(\begin{array}{cc} \cos\psi & -\sin\psi \\ \sin\psi &
      \cos\psi \end{array}\right)
\end{equation}
is a rotation matrix that describes the rotation of the car's fixed coordinate system 
against the inertial coordinate system. 

Under the assumption that the lateral velocity components at rear axle and
front axle vanish, we have that $r_R$, if transformed to the car's reference
system, only has a velocity component in the longitudinal direction of the
car, i.e. 
\begin{displaymath}
  \left(\begin{array}{c} v \\ 0 \end{array}\right) = S(\psi)^\top v_R \quad 
  \Longleftrightarrow \quad 
  v_R = S(\psi) \left(\begin{array}{c} v \\ 0 \end{array}\right) = 
  \left(\begin{array}{c} v \cos\psi \\ v \sin\psi \end{array}\right).
\end{displaymath}
This leads to the differential equations for the position $(x,y)$ of the
midpoint of the rear axle: 
\begin{displaymath}
  x'(t) = v(t) \cos \psi(t),\qquad y'(t) = v(t) \sin\psi(t).
\end{displaymath}
Under the assumption that the lateral velocity component at the front
axle vanishes, we have 
\begin{displaymath}
  0 = e_q^\top v_{F,b},
\end{displaymath}
where $e_q = (-\sin\delta,\cos\delta)^\top$ denotes the lateral direction of
the front axle and 
\begin{displaymath}
  v_{F,b} := S(\psi)^\top v_F = \left(\begin{array}{c} v \\ \ell \psi' \end{array}\right)
\end{displaymath}
denotes the representation of the velocity $v_F$ in the car's body fixed
reference system. The equation $0 = e_q^\top v_{F,b}$ yields the differential 
equation
\begin{displaymath}
  \psi'(t) = \frac{v(t)}{\ell} \tan\delta(t)
\end{displaymath}
for the yaw angle $\psi$. In summary, given the velocity $v(t)$ and the steering angle $\delta(t)$, 
the car's motion is given by the following system of differential equations: 
\begin{eqnarray}
  x'(t) & = & v(t) \cos \psi(t), \label{EQ:1a}\\
  y'(t) & = & v(t) \sin\psi(t), \label{EQ:1b}\\
  \psi'(t) & = & \frac{v(t)}{\ell} \tan\delta(t). \label{EQ:1c}
\end{eqnarray}

The steering angle and the velocity are typically bounded by 
$| \delta | \leq \delta_{max}$ and by $0 \leq v \leq v_{max}$ with given bounds $\delta_{max}$ and 
$v_{max}$, respectively. Moreover, one or more of the following modifications and restrictions 
can be used to yield a more realistic motion of the car: 
\begin{itemize}
\item
  Instead of controlling the velocity directly one often controls the acceleration
  by adding the differential equation
  \begin{equation}\label{EQ:2a}
    v'(t) = a(t)
  \end{equation}
  with bounds for the control $a$, i.e. $a_{min} \leq a \leq a_{max}$.
\item
  In order to model a certain delay in controlling the velocity $v$, one can 
  consider the following differential equation for $v$:
  \begin{equation}\label{EQ:2aa}
    v'(t) = \frac{v_{d}(t) - v(t)}{T}\qquad (T>0)
  \end{equation}
  Herein, $v_d(t)$ is the reference (=desired) velocity viewed as a control input 
  and $v$ is the actual velocity. The constant $T$ allows to influence the 
  response time, i.e. the delay of $v$. 
\item
  Instead of controlling the steering angle directly one often controls the 
  steering angle velocity by adding the differential equation
  \begin{equation}\label{EQ:2b}
    \delta'(t) = w(t)
  \end{equation}
  with bounds for the control $w$, i.e. $|w| \leq w_{max}$.
\end{itemize}

The car model derived in this section only provides a simple model, which is, however, very 
useful for path planning tasks. More detailed models can be found in \cite{Rill2011}.

\section{Trajectory Optimization and Path Planning}
\label{sec:trajectory}

\subsection{Geometric Path Planning by Shortest Paths}

A first approach towards the automatic path planning is based on purely geometric 
considerations, whereas the detailed dynamics of the vehicle are not taken into account in its full complexity. 
In fact, a shortest path is sought that connects two given points in a configuration space, e.g. vehicle positions, 
while taking into account a fixed number of obstacles. The result is a sequence of configuration points 
that need to be tracked by the vehicle. This simple, but robust and effective planning strategy 
turns out to be useful in the presence of complicated path constraints or 
many fixed obstacles that need to be avoided. 
The resulting path may serve as a reference path for inverse dynamics, a feedback tracking controller, 
or as an initial guess for more sophisticated optimization approaches from optimal control. 

In summary, the following steps need to be performed in order to compute a trajectory from an initial 
position to a terminal position, while taking obstacles into account: 
\begin{itemize}
\item[(1)]
  Solve a collision-free geometric shortest path problem using Dijkstra's algorithm. 
  The shortest path consists of a sequence of way-points leading from the initial 
  position to the terminal position.
\item[(2)]
  Interpolate the way-points by a cubic spline function (use a thinning algorithm if 
  necessary).
\item[(3)]
  Use inverse kinematics based on the car model or a feedback controller to track the spline 
  function. 
\end{itemize}
The above steps are discussed step by step. 

In a first attempt we focus on the two dimensional $(x,y)$-plane as the configuration space, 
which is often sufficient for autonomous ground-based vehicles. Note that the subsequent techniques can be easily 
extended to higher dimensions, e.g. in the context of robotics or flight path optimization. 

Consider a two dimensional configuration space 
\begin{displaymath}
  Q=\{ (x,y)^\top \in\R^2 \; |\; x_{min} \leq x\leq x_{max}, y_{min}\leq y\leq y_{max}\}
  = [x_{min},x_{max}]\times [y_{min},y_{max}],
\end{displaymath}
which corresponds to the feasible $(x,y)$-positions of the reference point on the vehicle. 
An equidistant discretization of the configuration space reads as
\begin{displaymath}
  Q_h = \{ (x_i,y_j)^\top \in \R^2\;|\; x_i = x_{min} + i h_x, \; 
  y_j = y_{min} + j h_x, i=0,\ldots,N_x, j=0,\ldots,N_y\} 
\end{displaymath}
with step-sizes $h_x = (x_{max}-x_{min})/N_x$ and $h_y = (y_{max}-y_{min})/N_y$ and given numbers 
$N_x,N_y\in\N$. Every grid point $q\in Q_h$ corresponds to the position of 
the vehicle's reference point. 

\begin{figure}
  \begin{center}
    \includegraphics[width=8cm]{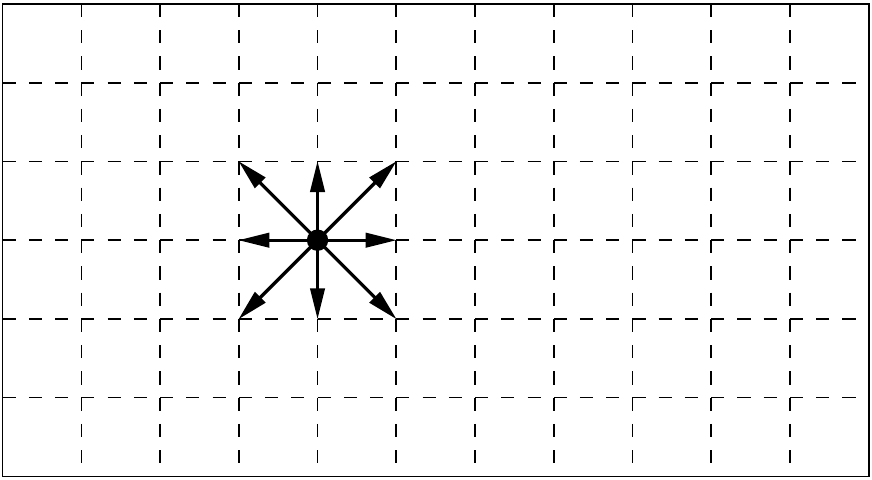}
    \begin{picture}(0,0)
      \put(-240,60){$y$}      
      \put(-117,-10){$x$}
    \end{picture}
  \end{center}
  \caption{Configuration space $Q_h$ and feasible transitions from a given grid point.}\label{Fig:ShortestPath}
\end{figure}
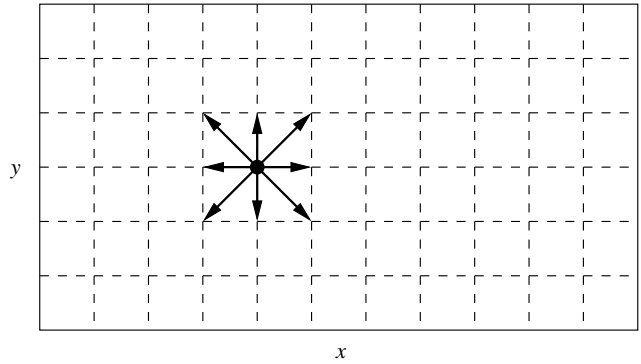

The geometric motion of the vehicle can be modeled in a 
very simplified way by transitions from a given grid point to neighboring 
grid points (called feasible transitions), see Figure~\ref{Fig:ShortestPath}. 

To this end, the discrete configuration space $Q_h$ together with the feasible transitions, defined by 
the discrete control set $\mathcal{U}_h$, define a directed graph $G = (Q_h,E_h)$ with nodes $Q_h$ and 
edges $E_h \subset Q_h\times Q_h$ such that 
\begin{displaymath}
  ( \bar q, q )  \in E_h \qquad \Longleftrightarrow \qquad | x - \bar x | \leq h_x \mbox{ and } |y-\bar x| \leq h_y ,
\end{displaymath}
where $\bar q = (\bar x,\bar y)^\top\in Q_h$ and $q = (x,y)^\top\in Q_h$. Note that the set of edges can be further 
reduced if only some of the transitions between grid points are permitted. 

To each edge $(\bar q,q)\in E_h$ of the graph we may assign the cost $c(\bar q,q) = \| q - \bar q\|_2$ assuming that 
the time to move from node $\bar q$ to node $q$ is proportional to its Euclidean distance. More general 
costs may be assigned as well as long as they are non-negative. In particular we have to deal with infeasible nodes 
owing to collisions with obstacles. 
In order to avoid collisions, infeasible nodes can be eliminated from the set $Q_h$ beforehand or an
infinite cost $c(\bar q,q)$ can be assigned whenever $\bar q$ or $q$ are infeasible nodes. Whether a collision 
occurs can be checked with the technique of Section~\ref{Sec:CollisionDetection}. 

Given the graph $G$ and the non-negative cost function $c: E_h \rightarrow \R$, the task is to find a shortest path in $G$ 
leading from an initial node $q_0 \in Q_h$ to a terminal node $q_T\in Q_h$. This can be achieved by Dijkstra's 
shortest path algorithm or versions of it, compare \cite{Pap98,Kor08}. The algorithm, which exploits the 
dynamic programming principle, has a complexity of $\mathcal{O}(n^2)$, where $n$ denotes the number of nodes in $G$. 
An efficient implementation of Dijkstra's algorithm uses a priority queue, see \cite{Cormen2009}.

After a shortest path has been found, it can be further tuned in a post-processing process. 
The post-processing includes an optional thinning of the path as described in 
\cite{Nachtigal2015}, a spline interpolation, and the application of inverse kinematics or 
a tracking controller to follow the spline with the actual car model. A tracking controller 
based on the car model in (\ref{EQ:1a})-(\ref{EQ:1c}) is designed in Section~\ref{sec:Tracking}.

\begin{figure}
  \begin{center}
    \includegraphics[height=\textwidth,angle=-90]{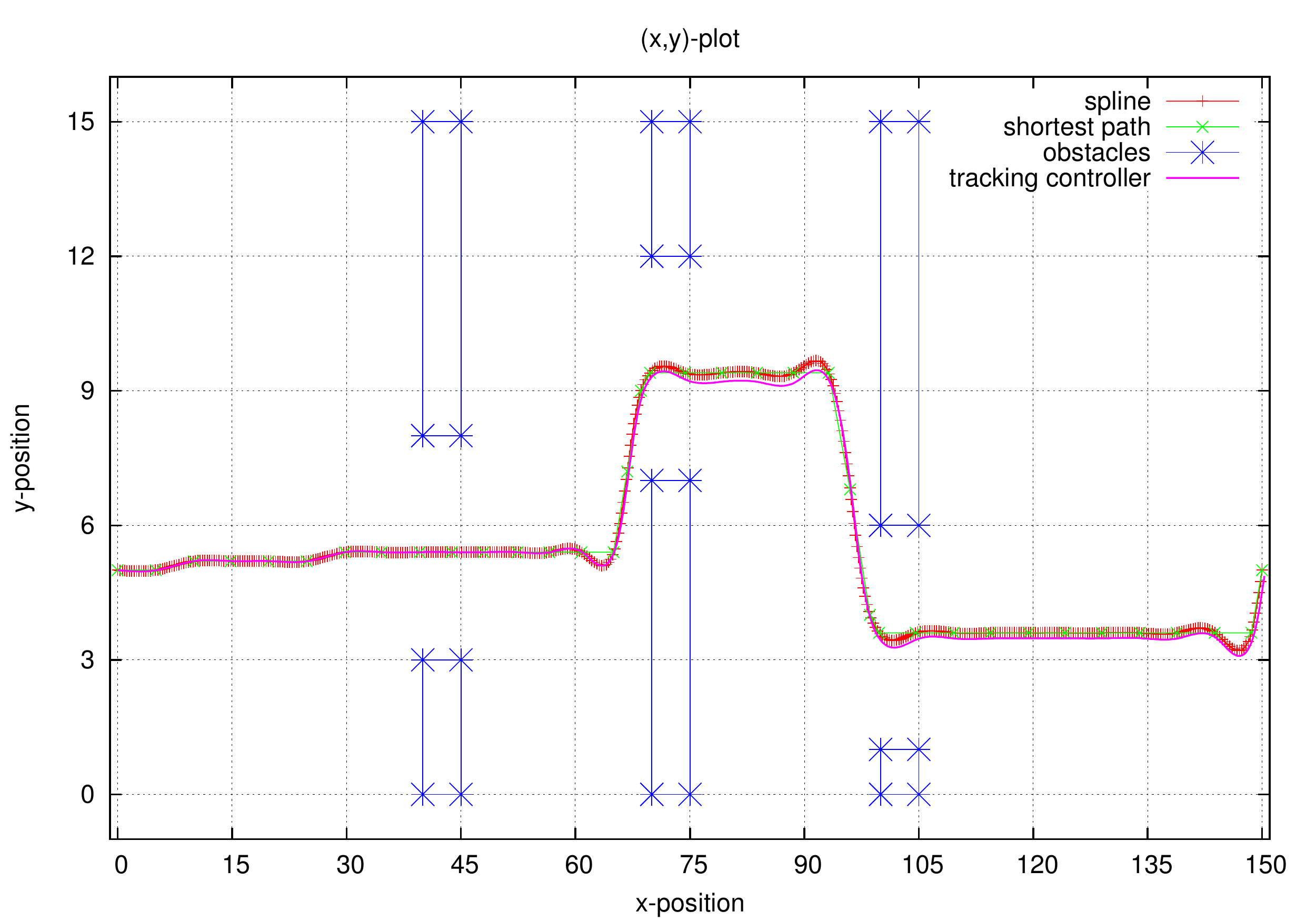}
  \end{center}
  \caption{A car's trajectory in the $(x,y)$-plane in the presence of obstacles  
    starting at initial position $(0,5)$ and ending at position $(150,0)$: Shortest path, 
    spline interpolation and tracked curve. }
  \label{Fig:DijkstraLaneChange}
\end{figure}

Figure~\ref{Fig:DijkstraLaneChange} shows the result of the geometric shortest path 
planning approach for a double lane change maneuver. The trajectory is typically not
optimal in the sense of optimal control, but may be acceptable for practical purposes or 
it may serve as an initial guess for an optimal control problem.  

\begin{figure}
  \begin{center}
    \includegraphics[height=\textwidth,angle=-90]{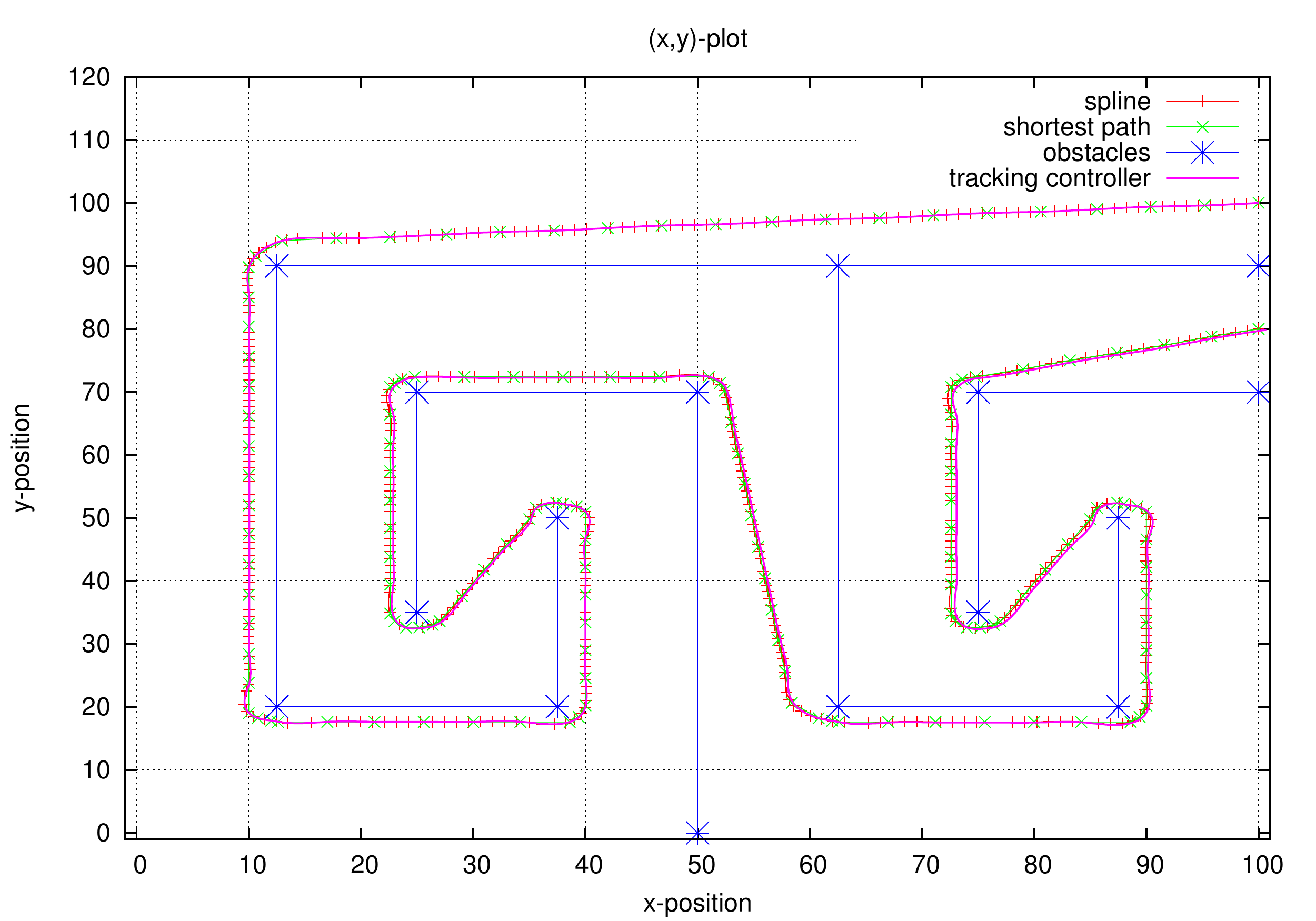}
  \end{center}
  \caption{A car's trajectory in the $(x,y)$-plane in the presence of obstacles  
    starting at initial position $(100,100)$ and ending at position $(100,80)$: 
    Shortest path, spline interpolation and tracked curve. }
  \label{Fig:DijkstraScenario2}
\end{figure}

A more complicated track is depicted in Figure~\ref{Fig:DijkstraScenario2}. 

Although the geometric shortest path problem is robust and can handle fixed obstacles, 
moving obstacles (at the cost of an additional time state variable), and 
complicated obstacle geometries very well, it suffers from the disadvantage that, e.g. the 
steering angle $\delta$ cannot be constrained. For instance, the course in 
Figure~\ref{Fig:DijkstraScenario2} requires a maximum steering angle of 73 degrees, which is 
beyond the technical limit of a car. To avoid this remedy, it is possible to compute a cubic 
spline approximation in combination with inverse kinematics subject to constraints for the 
steering angle. Alternatively, one can modify the shortest path approach as follows: 
In the previous model, a purely geometric motion was considered that can be interpreted as 
controlling the motion in $x$- and $y$-direction independently from each other. 
This might be reasonable for omnidirectional robots, but it does not reflect the actual motion 
capabilities of a car with a steering device very well. To this end it is more realistic, 
although computationally more expensive, to work with a discretization of the dynamics 
(\ref{EQ:1a})-(\ref{EQ:1c}) and the control set 
$\mathcal{U}:=[v_{min},v_{max}]\times [\delta_{min},\delta_{max}]$. This leads to a three 
dimensional configuration space $Q$ with points of type $q = (x,y,\psi)^\top$ subject to suitable 
bounds on the components. The discretized configuration space $Q_h$ is constructed as follows,
where we use the explicit Euler method for simplicity: Let $\bar q=(\bar x,\bar y,\bar \psi)^\top\in Q_h$ 
be arbitrary. Then $Q_h$ contains all the points $q=(x,y,\psi)^\top \in Q$ of type  
\begin{displaymath}
  q = \left(\begin{array}{c} x \\ y \\ \psi \end{array}\right) 
    = \left(\begin{array}{c}
      \bar x + h v \cos \bar \psi \\
      \bar y + h v \sin \bar \psi \\
      \bar\psi + h \frac{v}{\ell} \tan\delta
    \end{array}\right) =: F(\bar q,u)
\end{displaymath}
with $u = (v,\delta)^\top \in \mathcal{U}_h$, where 
\begin{displaymath}
  \mathcal{U}_h := \{ (v_i,\delta_j)^\top \in \R^2\;|\; v_i = v_{min} + i h_v, \; 
  \delta_j = \delta_{min} + j h_\delta, i=0,\ldots,N_v, j=0,\ldots,N_\delta\}
\end{displaymath}
is a discrete approximation to the set $\mathcal{U}$ with step-sizes $h_v = (v_{max}-v_{min})/N_v$ and 
$h_\delta = (\delta_{max}-\delta_{min})/N_\delta$.

Again, the discrete configuration space $Q_h$ together with the feasible transitions, defined by 
the discrete control set $\mathcal{U}_h$, define a directed graph $G = (Q_h,E_h)$ with nodes $Q_h$ and 
edges $E_h \subset Q_h\times Q_h$ such that 
\begin{displaymath}
  ( \bar q, q )  \in E_h \qquad \Longleftrightarrow \qquad \exists u\in \mathcal{U}_h : q = F(\bar q,u) . 
\end{displaymath}
The graph is considerably larger than the previous one and the computational effort for solving 
the corresponding shortest path problem increases accordingly. Please note that time dependent 
obstacle motions can be incorporated at the dispense of an additional configuration variable 
which corresponds to the time. The shortest path approach essentially coincides with the 
dynamic programming approach, compare \cite{Bertsekas2012}.

\subsection{Collision Detection}
\label{Sec:CollisionDetection}

Collision detection is an important issue in robotics and autonomous systems and various approaches 
based on distance functions of convex bodies have been developed, see \cite{Gilbert1985,Johnson1985,Gilbert1988}. A smoothed distance measure was constructed in
\cite{Escande2014} and can be used in gradient type optimization algorithms. 

Often it is sufficient to approximate the vehicle with center $(x_c,y_c)$ and obstacles with 
centers $(x_i,y_i)$ by circles of radius $r$ and $r_i$, $i=1,\ldots,M$, respectively, 
and to impose state constraints of type 
\begin{displaymath}
  (x_c - x_i)^2 + (y_c - y_i)^2 \geq (r+r_i)^2,\qquad i=1,\ldots,M,
\end{displaymath}
to prevent collisions. Herein, the center of the vehicle based on the configuration in 
Figure~\ref{Fig:1} is given by
\begin{equation}\label{EQ:CarCenter}
  r_c := \left(\begin{array}{c} x_c \\ y_c \end{array}\right) = \left(\begin{array}{c} x \\ y \end{array}\right)
  + S(\psi) \left(\begin{array}{c} \ell/2 \\ 0 \end{array}\right). 
\end{equation}
Collision detection, taking the detailed shape of the vehicle 
and the obstacles into account, is much more involved, since it is not straightforward to compute the distance function 
for potentially non-convex bodies. We follow a technique in \cite{Landry2012} and assume that the shape of the vehicle of length 
$\ell$ and width $w$ in the two dimensional plane is given by a rectangle (w.r.t. to the vehicle's 
coordinate system)
\begin{displaymath}
  R = \{ z\in \R^2 \; |\; A z\leq b\}, \qquad 
  A = \left(\begin{array}{rr} 1 & 0 \\ -1 & 0 \\ 0 & 1 \\ 0 & -1 \end{array}\right),\ 
  b = \left(\begin{array}{r} \ell/2 \\ -\ell/2 \\ w/2 \\ -w/2 \end{array}\right).
\end{displaymath}
The rectangle moves along with the vehicle and its location at time $t$ is given by 
\begin{displaymath}
  R(t) = S(\psi(t)) R + r_c(t) = \{ z\in\R^2 \; |\; A S(\psi(t))^\top z \leq b + A S(\psi(t))^\top r_c(t)\}
\end{displaymath}
with $S$ from (\ref{EQ:Drehmatrix}) and $r_c$ from (\ref{EQ:CarCenter}). Let the $M$ obstacles at time $t$ 
be given by the union of convex polyhedra
\begin{displaymath}
  Q_i(t) := \bigcup_{j=1}^{M_i} Q^{(i,j)}(t) \qquad \text{with} \qquad Q^{(i,j)}(t)=\{ y\in \R^2 \; |\; C^{(i,j)}(t) y\leq d^{(i,j)}(t)\},
\end{displaymath}
where $M_i$ is the number of polyhedra in obstacle $Q_i$ and for $j=1,\ldots,M_i$, the matrix 
$C^{(i,j)}(t)\in \R^{q_{i,j}\times 2}$ and the vector $d^{(i,j)}(t)\in \R^{q_{i,j}}$ define the convex parts of the 
i-th obstacle at time $t$. Herein, $q_{i,j}$ is the number of facets in $Q^{(i,j)}$. 

The vehicle and the obstacles do not collide at time $t$ if and only if  
\begin{displaymath}
  R(t) \cap Q^{(i,j)}(t)  \,=\,\emptyset \qquad \forall j=1,\ldots,M_i, i=1,\ldots,M.
\end{displaymath}
This is equivalent to the infeasibility of the linear system
\begin{equation}\label{EQ:intersection}
  \left(\begin{array}{c}
    A(t)  \\ C^{(i,j)}(t) 
  \end{array}\right)\, z\,\leq\,\left(\begin{array}{c}
    b(t) \\ d^{(i,j)}(t)
  \end{array}\right),\quad \forall j=1,\ldots,M_i, i=1,\ldots,M,
\end{equation}
where $A(t) := A S(\psi(t))^\top$ and $b(t) := b + A S(\psi(t))^\top r_c(t)$. According to the 
Lemma of Gale the system (\ref{EQ:intersection}) has no solution at time $t$ if and only if there exist 
 vectors $ w^{(i,j)}(t)\in\R^{4+q_{i,j}}$ such that 
\begin{displaymath}
  w^{(i,j)}(t)\,\geq\,0,\quad
  \left(\begin{array}{c}
    A(t) \\ C^{(i,j)}(t)
  \end{array}\right)^\top w^{(i,j)}(t)\,=\,0 \quad \text{and}\quad 
  \left(\begin{array}{c}
    b(t) \\ d^{(i,j)}(t)
  \end{array}\right)^\top w^{(i,j)}(t) < 0.
\end{displaymath}
Note that a scaled vector $\lambda w^{(i,j)}(t)$ with $\lambda>0$ satisfies the conditions as well, if $ w^{(i,j)}(t)$ does so. 
Hence, we may bound the length of the components by one and impose the additional constraint $w^{(i,j)}(t) \leq e$,
where $e = (1,\ldots,1)^\top$ denotes the vector of all ones of appropriate dimension.

This condition can be checked by solving the following linear program for all 
$j=1,\ldots,M_i$, $i=1,\ldots,M$, and all $t$: 

\medskip
{\em Minimize 
\begin{displaymath}
  \left(\begin{array}{c}
    b(t) \\ d^{(i,j)}(t)
  \end{array}\right)^\top w^{(i,j)}(t)
\end{displaymath}
\indent subject to the constraints
\begin{displaymath}
  0\,\leq \, w^{(i,j)}(t) \,\leq\,e,\qquad
  \left(\begin{array}{c}
    A(t)\\C^{(i,j)}(t)
  \end{array}\right)^\top w^{(i,j)}(t)\,=\,0.
\end{displaymath}}

\noindent
Note that the feasible sets of the linear programs are non-empty and compact and thus an 
optimal solution exists.  

A collision does not occur, if the value function 
\begin{displaymath}
  \zeta^{(i,j)}(t) := \min\left\{  \left(\begin{array}{c}
    b(t)^\top \\ d^{(i,j)}(t)
  \end{array}\right)^\top w \;\Bigg|\;  0\,\leq w\,\leq\,e,\;
  \left(\begin{array}{c}
    A(t) \\C^{(i,j)}(t)
  \end{array}\right)^\top w=0\right\}
\end{displaymath}
is negative for all combinations $(i,j)$ and all $t$. Hence, collisions are avoided by imposing 
the non-linear and non-differentiable constraint
\begin{displaymath}
  \sup_{t\in [0,t_f],(i,j)} \zeta^{(i,j)}(t) \leq -\varepsilon
\end{displaymath}
for some $\varepsilon>0$ sufficiently small. Note that $\zeta^{(i,j)}$ implicitly depends on the vehicle's state $(x,y,\psi)$.

As an alternative to the above linear programming approach, the Gilbert-Johnson-Keehrti 
algorithm from \cite{Gilbert1988} is frequently used in computer graphics and robotics for 
real-time collision detection. It uses Minkowski sums and convex hulls to compute the 
signed distance function between two polyhedral objects.

\subsection{Optimal Drivers by Optimal Control}
\label{sec:OCP}

A virtual ``optimal'' driver can be modeled by means of a suitable optimal control problem, which fits 
into the following general class of parametric optimal control problems OCP($p$): 

\medskip
{\em Minimize 
\begin{equation} \label{EQ:OCP:OBJ}
  \varphi(z(t_f),t_f,p) + \int_{t_0}^{t_f} f_0(z(t),u(t),p) dt
\end{equation}
\indent subject to the constraints
\begin{align}
  z'(t) & = f(z(t),u(t),p) && \mbox{a.e. in } [t_0,t_f], \label{EQ:OCP:1}\\
  g(t,z(t),p) & \leq 0 && \mbox{in } [t_0,t_f], \label{EQ:OCP:2}\\
  \psi(z(t_0),z(t_f),p) & = 0, \label{EQ:OCP:3} \\
  u(t) & \in \mathcal{U} & & \mbox{a.e. in } [t_0,t_f]. \label{EQ:OCP:4}
\end{align}}

Herein, the objective function (\ref{EQ:OCP:OBJ}) typically consists of a linear combination 
of final time, steering effort and fuel consumption. The car model defines the differential equation 
in (\ref{EQ:OCP:1}) for the state $z$ and control $u$, while road boundaries and stationary or 
moving obstacles lead to state constraints of type (\ref{EQ:OCP:2}). Initial and terminal states 
 of the 
vehicle at initial time $t_0$ and terminal time $t_f$ are limited by the constraint in (\ref{EQ:OCP:3}). Finally, the control vector $u$ is restricted 
to the control set $\mathcal{U}$ in (\ref{EQ:OCP:4}). The problem formulation may depend on a parameter vector $p$ that 
can be used to model perturbations or uncertainties, which enter the optimal control problem as 
parameters. 

Given a nominal parameter $p^*$, various approaches exist to solve OCP($p^*$) numerically. The indirect 
solution approach exploits first order necessary optimality conditions, see \cite{Iof79}. The function 
space approach applies optimization procedures, i.e. gradient type methods in the function space 
setting of OCP($p^*$), compare~\cite{Pol73}. For highly nonlinear problems with complicated state and 
control constraints, direct discretization methods, see \cite{Boc84,Bet01,Gerdts2012}, are often preferred 
owing to their flexibility and robustness. In the following example, we use an optimal control problem to model a
parking maneuver of a car and solve the optimal control problem by the direct shooting method OCPID-DAE1, 
see \cite{OCPIDDAE1}.

\begin{example}[Parking maneuver]
  \label{ex:Parken}
  The task is to park a car in a parking space next to the car on the road. 
  The car's dynamics are given by (\ref{EQ:1a})-(\ref{EQ:1c}), (\ref{EQ:2a}) and (\ref{EQ:2b}) with $\ell=2.7$ and width $b=1.8$. The state vector is given by $z=(x,y,\psi,v,\delta)^\top$ and the 
control vector by $u=(w,a)^\top$. 

  The steering angle velocity $w$ and 
  the acceleration $a$ are restricted by the control constraints
  \begin{equation}\label{EQ:Parken:1}
    w(t) \in [-0.5,0.5],\qquad a(t)\in [-0.5,0.5].
  \end{equation} 
  The steering angle $\delta$ is bounded by the state constraints 
  \begin{equation}\label{EQ:Parken:2}
    -\frac{\pi}{6} \leq \delta(t) \leq \frac{\pi}{6}.
  \end{equation}
  The initial state of the car is given by 
  \begin{equation}\label{EQ:Parken:3}
    x(0) = 2.5,\quad y(0) = 1.5,\quad \psi(0) = v(0) = \delta(0) = 0
  \end{equation}
  and the terminal state by 
  \begin{equation}\label{EQ:Parken:4}
    x(t_f) = -1.25,\quad y(t_f) = -1.5,\quad \psi(t_f) = v(t_f) = \delta(t_f) = 0.
  \end{equation}
  Moreover, the parking lot, which is located to the right of the car, is defined by the state constraints 
  \begin{equation}\label{EQ:Parken:5}
    y_r(t) \geq \eta(x_r(t)) \quad \mbox{and}\quad y_f(t) \geq \eta(x_f(t)), 
  \end{equation}
  where 
  \begin{displaymath}
    \left(\begin{array}{c} x_r(t) \\ y_r(t) \end{array}\right) = \left(\begin{array}{c} x(t) \\ y(t) \end{array}\right)
    + S(\psi(t)) \left(\begin{array}{c} 0 \\ b/2 \end{array}\right)\mbox{ and }
    \left(\begin{array}{c} x_f(t) \\ y_f(t) \end{array}\right) = \left(\begin{array}{c} x(t) \\ y(t) \end{array}\right)
    + S(\psi(t)) \left(\begin{array}{c} \ell \\ b/2 \end{array}\right)
  \end{displaymath}
  denote the positions of the right front and right rear wheel centers, respectively, $S$ is 
  the rotation matrix in (\ref{EQ:Drehmatrix}) and $\eta$ is the piecewise defined and continuously differentiable function
  \begin{displaymath}
    \eta(x) := \left\{\begin{array}{ll}
    0,  & \mbox{if } |x| \geq 2.5,\\
    -3, & \mbox{if } |x| \leq 2.4, \\
    -900 (|x| - 2.5)^2 - 6000 (|x|-2.5)^3, & \mbox{if } 2.4 < |x| < 2.5.
    \end{array}\right. 
  \end{displaymath}
  
  Summarizing, the optimal control problem aims at minimizing a linear combination of final time and steering effort, that is 
  \begin{displaymath}
    t_f + \int_0^{t_f} w(t)^2 dt, 
  \end{displaymath}
  subject to the constraints (\ref{EQ:1a})-(\ref{EQ:1c}), (\ref{EQ:2a}), (\ref{EQ:2b}), (\ref{EQ:Parken:1})-(\ref{EQ:Parken:5}).

  Figures~\ref{Fig:Parken:1a}-\ref{Fig:Parken:1b} show the result of the direct shooting method 
  OCPID-DAE1, see \cite{OCPIDDAE1}, with $N=101$ grid points and final time $t_f\approx 15.355$. 
  Please note the 3 phases of the parking maneuver in Figure~\ref{Fig:Parken:1a}. 

  \begin{figure}[h]
    \begin{center}
      \includegraphics[height=0.8\textwidth,angle=-90]{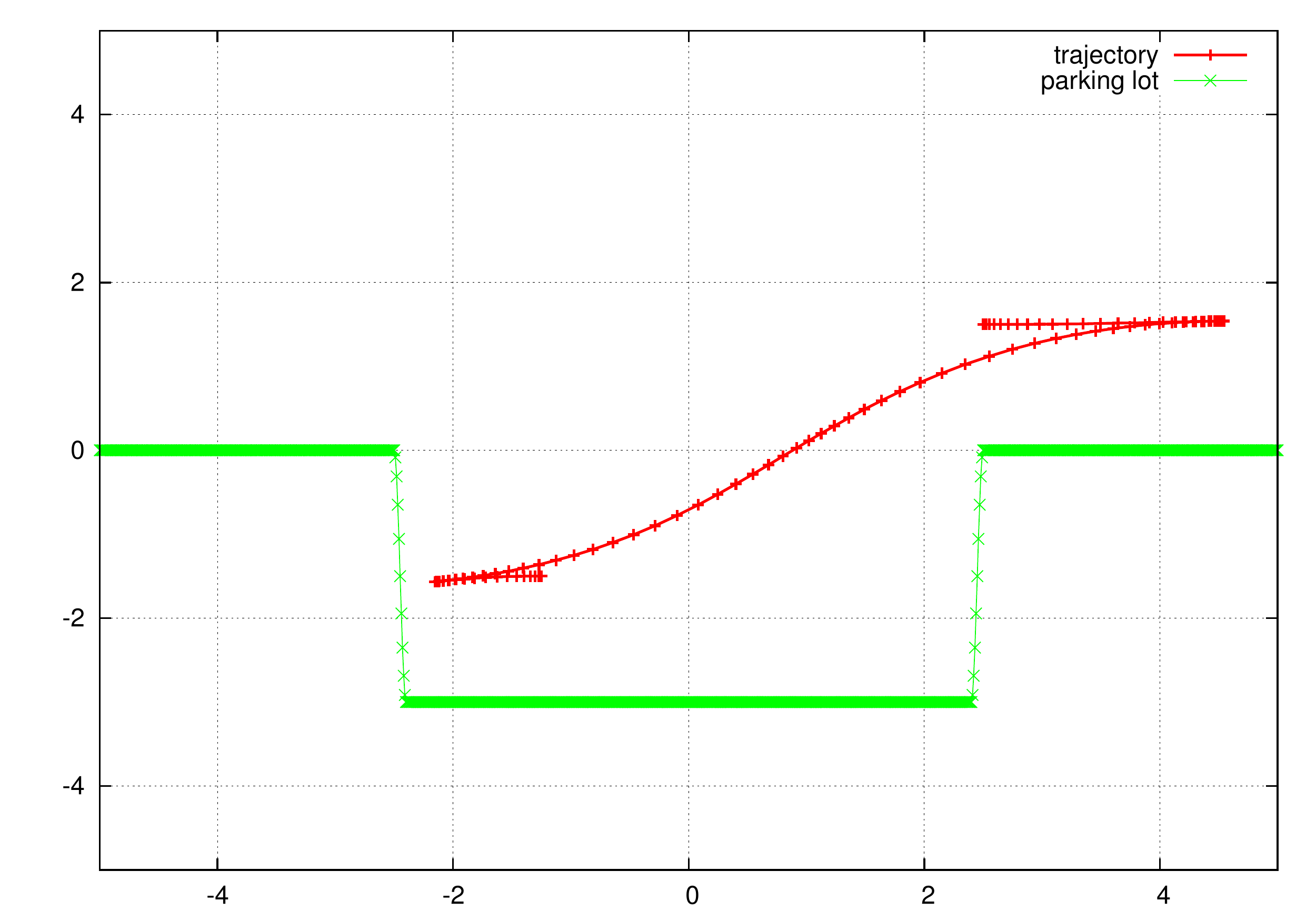}    
    \end{center}
    \caption{Trajectory $(x,y)$ for optimal parking maneuver.}
    \label{Fig:Parken:1a}
  \end{figure}

  \begin{figure}[h]
    \begin{center}
      \includegraphics[height=0.4\textwidth,angle=-90]{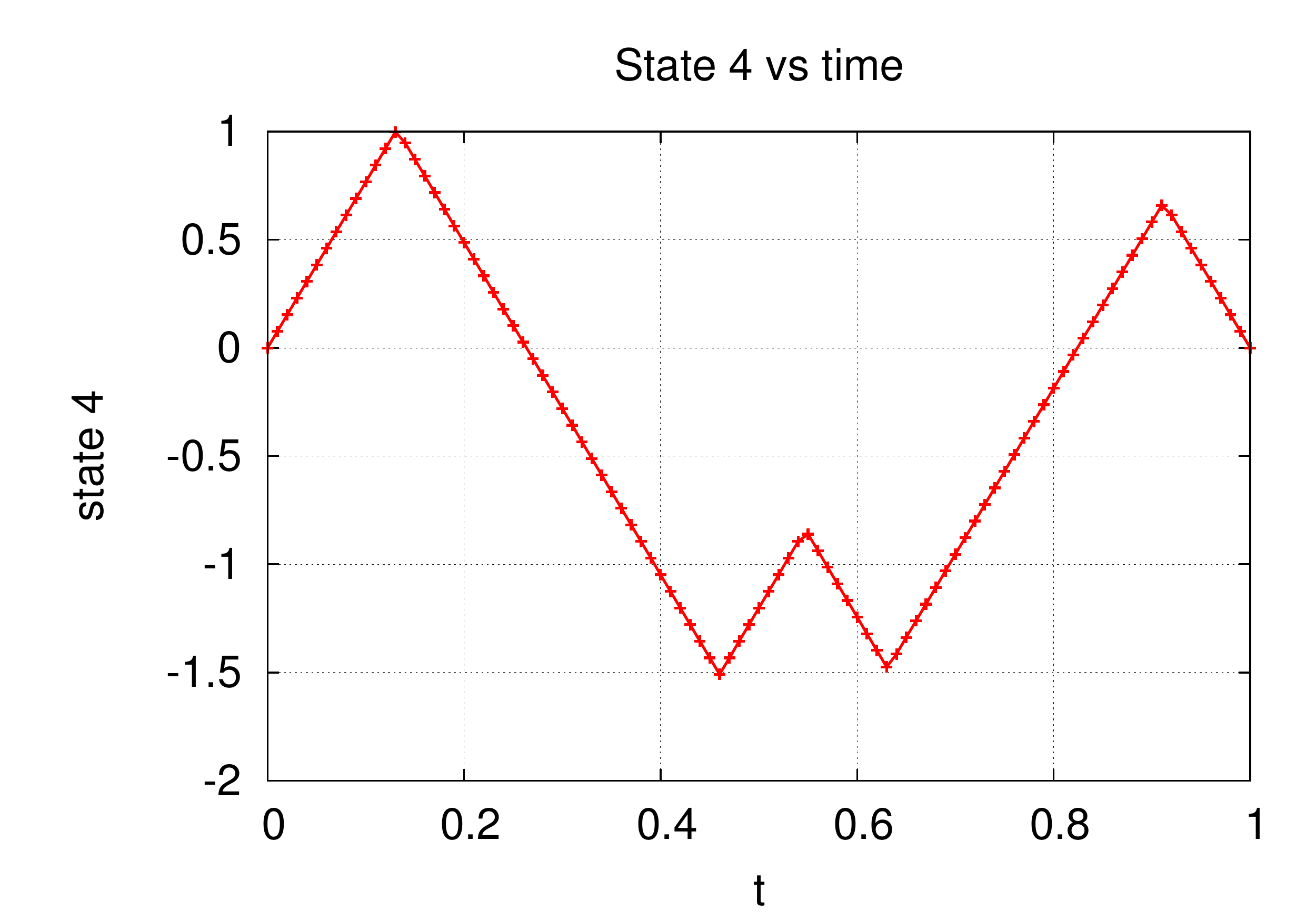}      
      \includegraphics[height=0.4\textwidth,angle=-90]{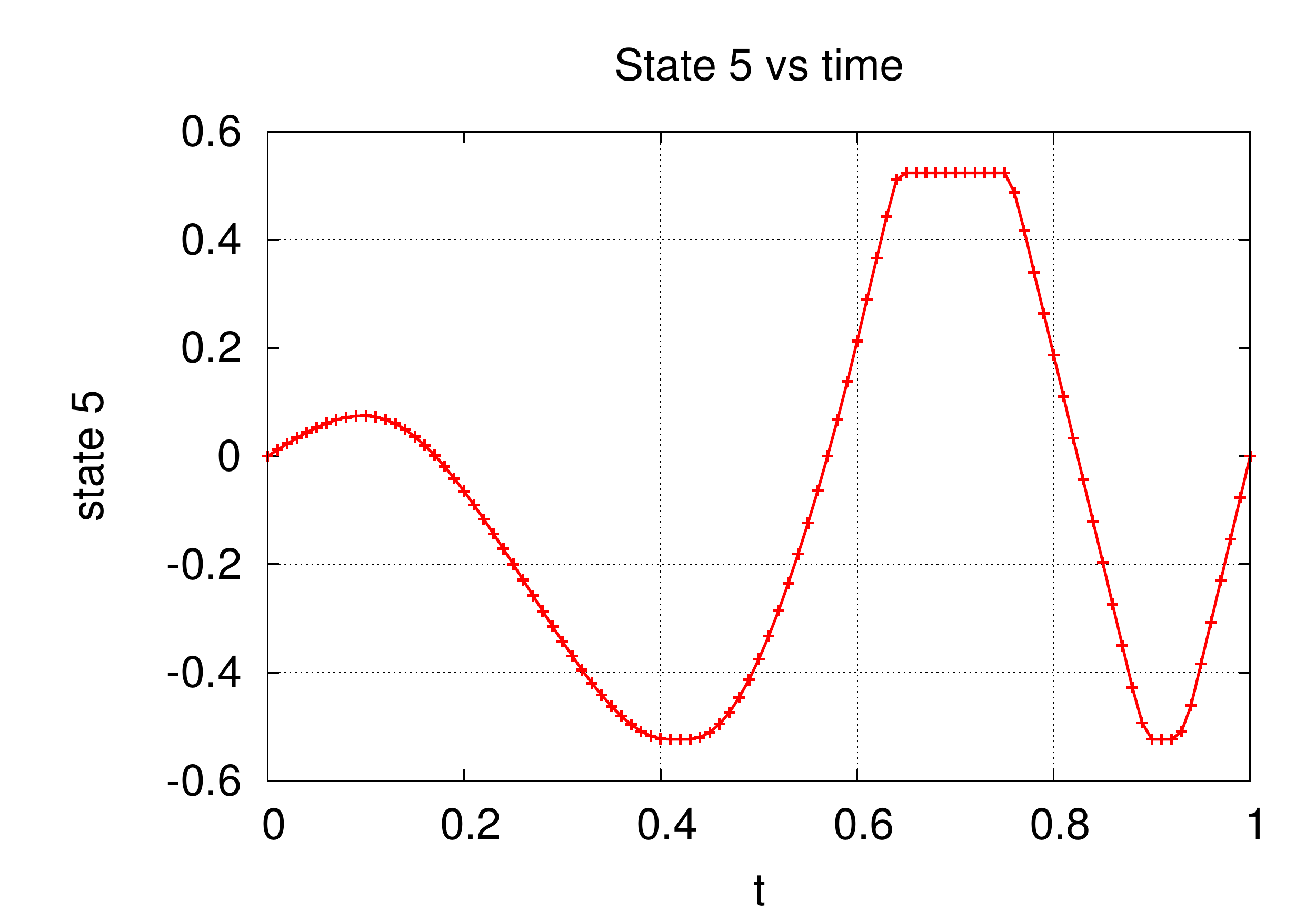}
      
      \includegraphics[height=0.4\textwidth,angle=-90]{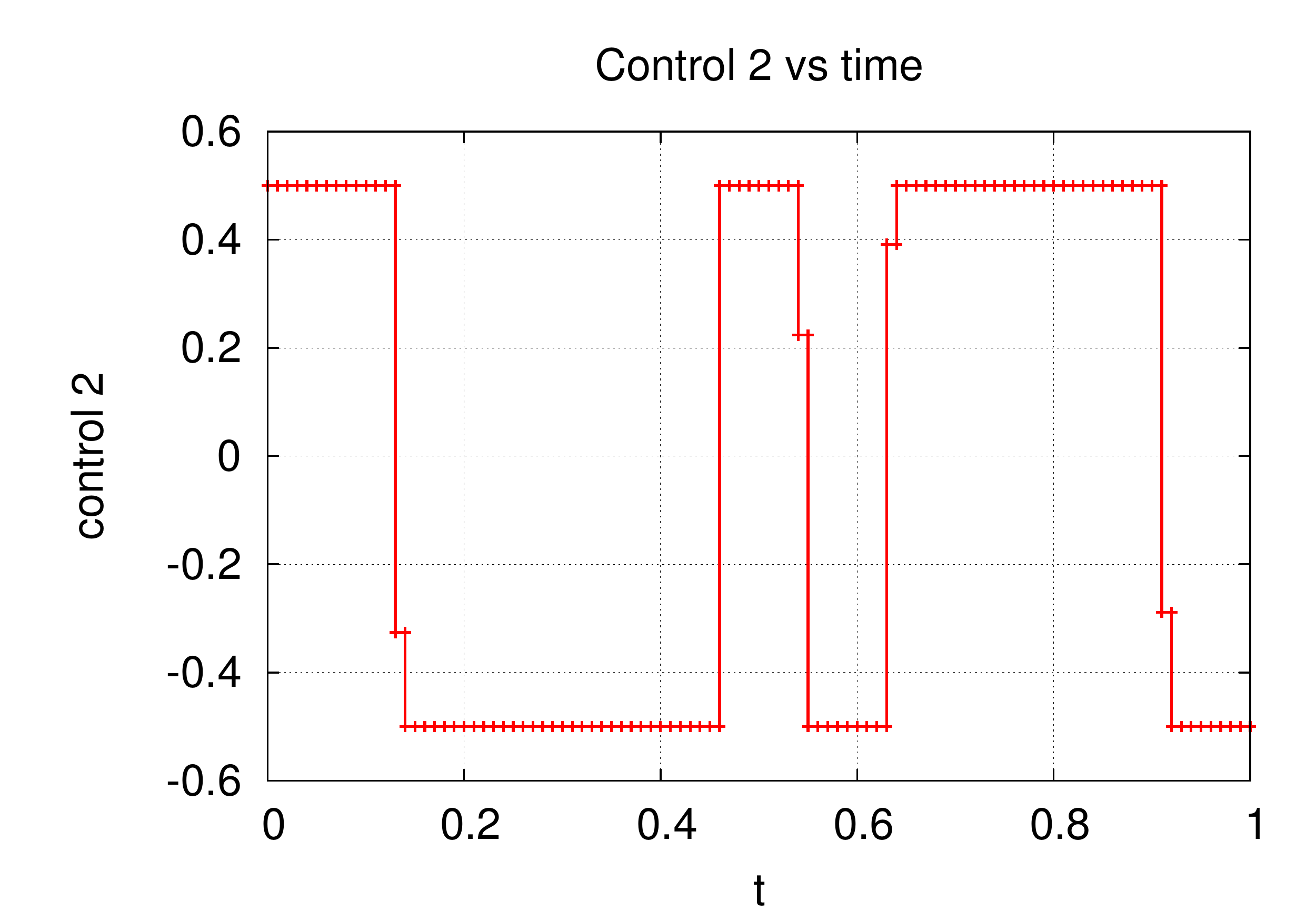}      
      \includegraphics[height=0.4\textwidth,angle=-90]{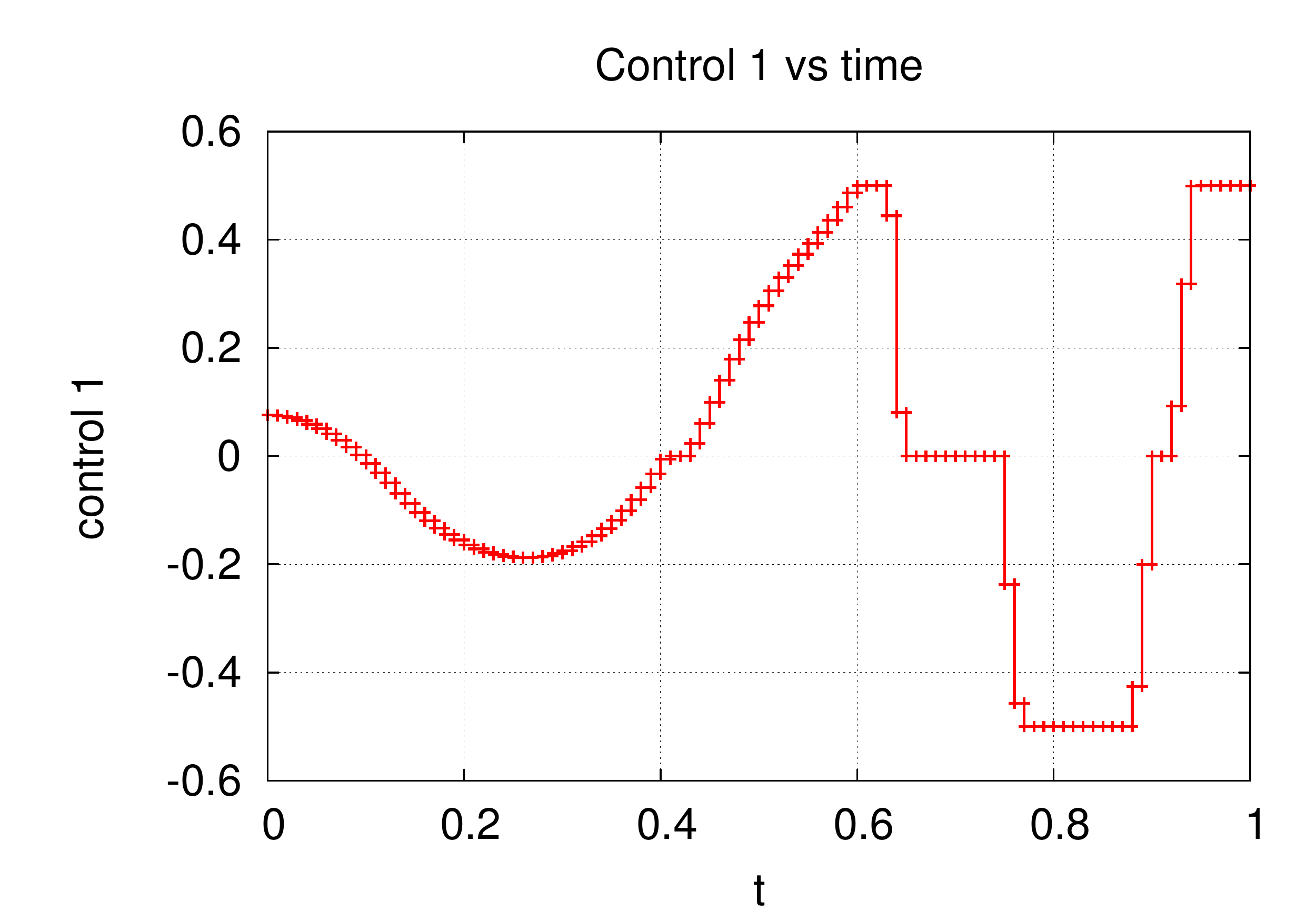}
    \end{center}

    \caption{Velocity $v$ (state, top left), steering angle $\delta$ (state, top right), 
      acceleration $a$ (control, bottom left), and steering angle velocity $w$ 
      (control, bottom right) for optimal parking maneuver.}
    \label{Fig:Parken:1b}
  \end{figure}
  
  Figure~\ref{Fig:Parken:2} illustrates the motion of the car with some snapshots. 

  \begin{figure}[h]
    \begin{center}
      \includegraphics[height=0.2\textwidth,angle=-90]{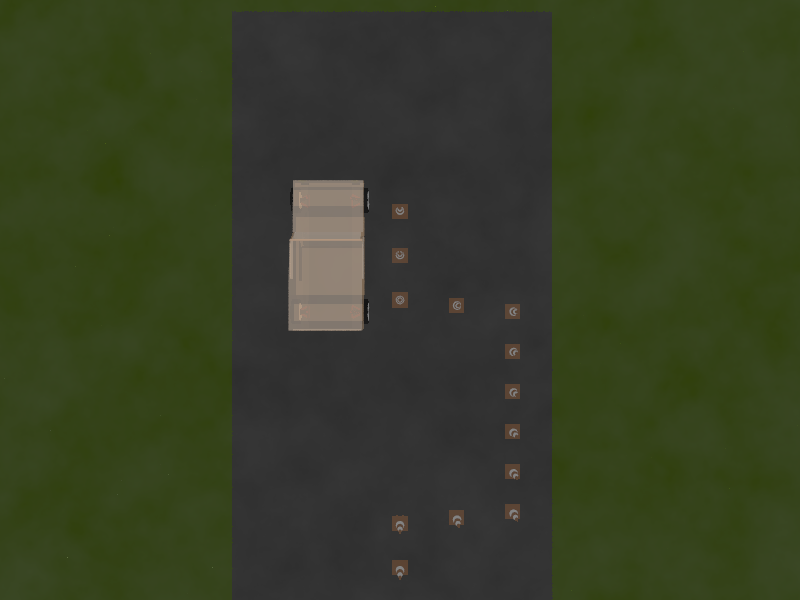}
      \includegraphics[height=0.2\textwidth,angle=-90]{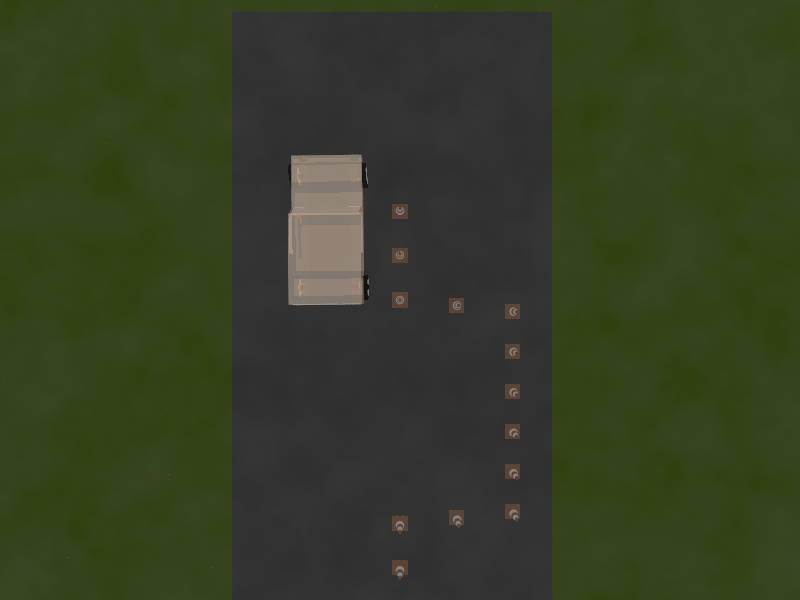}
      \includegraphics[height=0.2\textwidth,angle=-90]{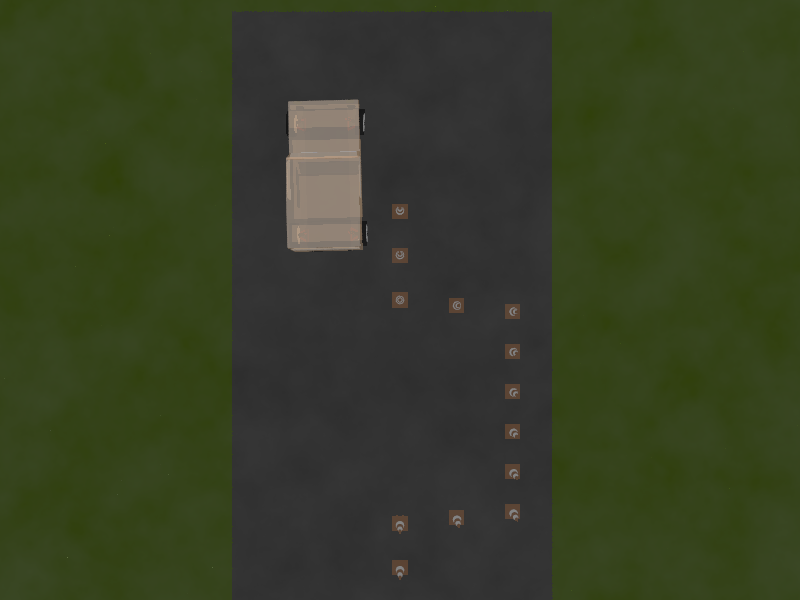}

      \includegraphics[height=0.2\textwidth,angle=-90]{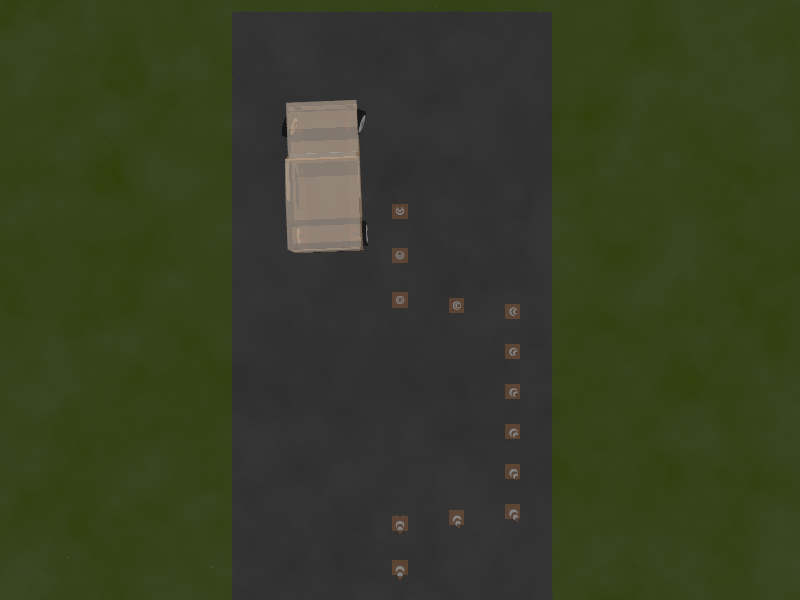}
      \includegraphics[height=0.2\textwidth,angle=-90]{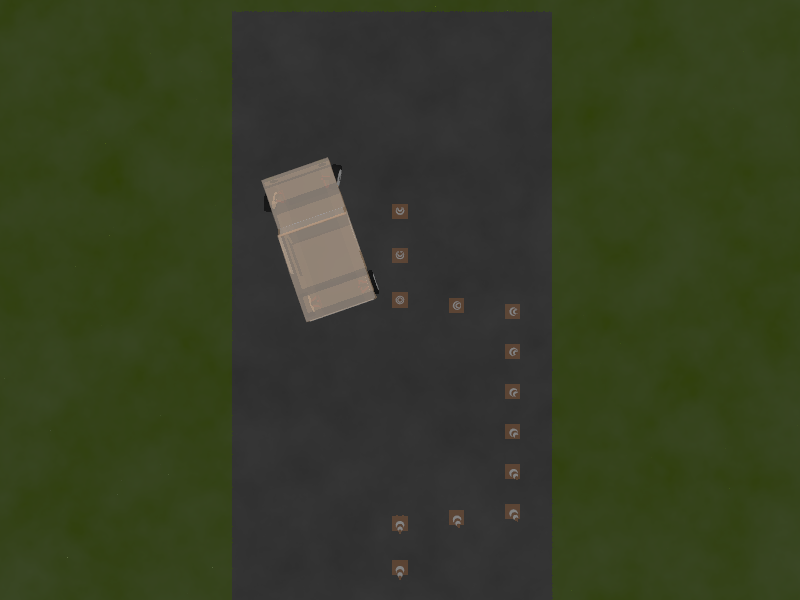}
      \includegraphics[height=0.2\textwidth,angle=-90]{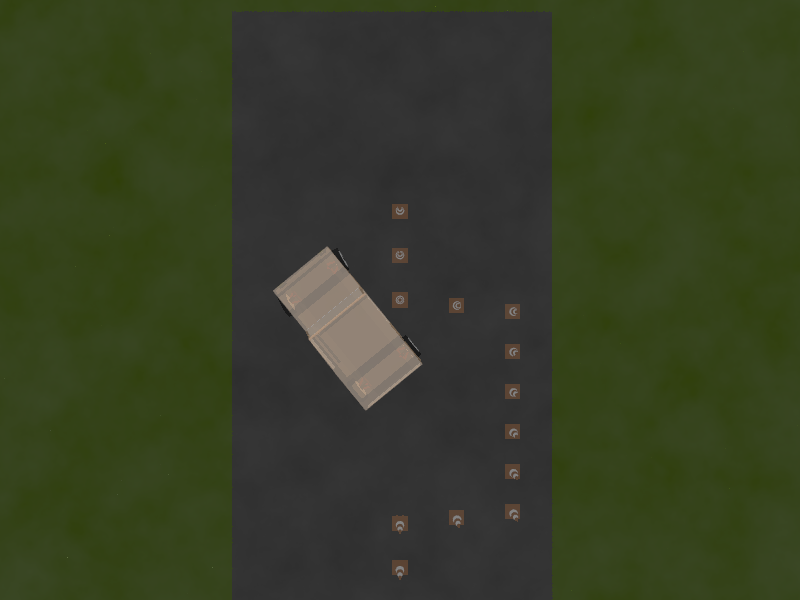}

      \includegraphics[height=0.2\textwidth,angle=-90]{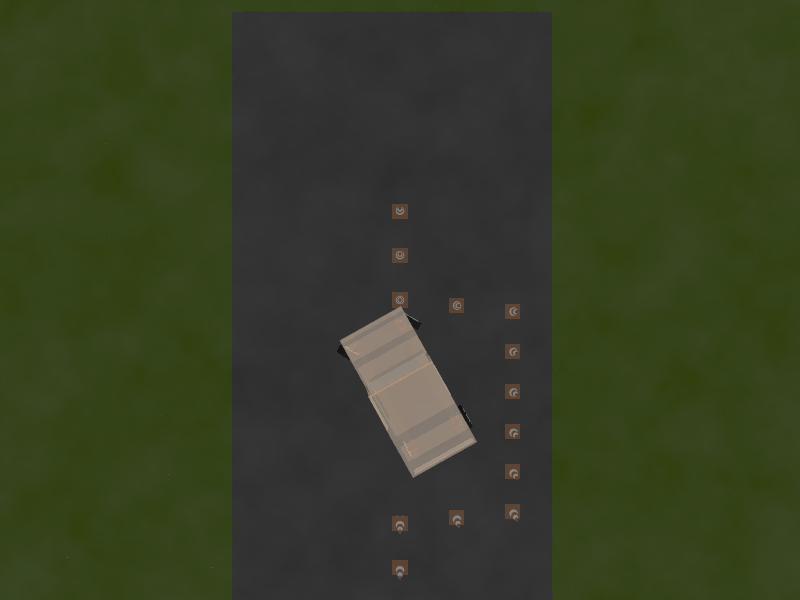}
      \includegraphics[height=0.2\textwidth,angle=-90]{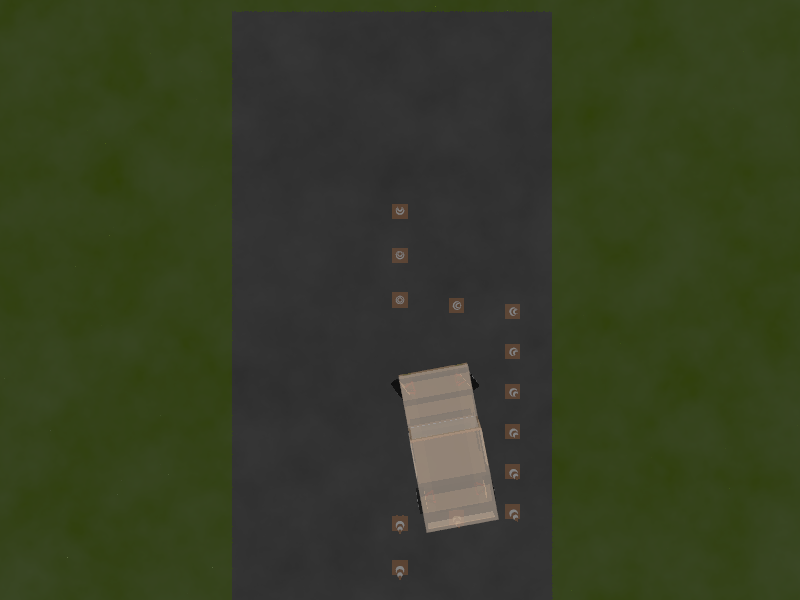}     
      \includegraphics[height=0.2\textwidth,angle=-90]{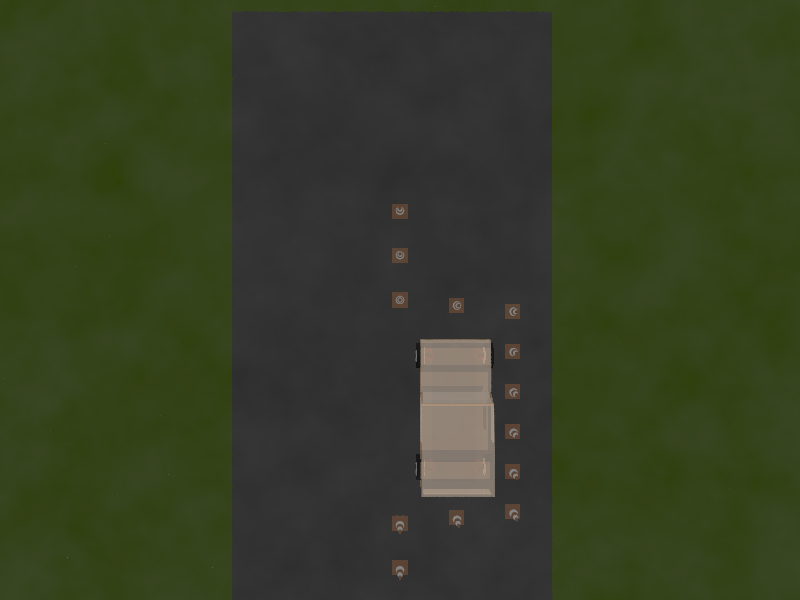}
    \end{center}

    \caption{Snapshots of the car's parking maneuver.}
    \label{Fig:Parken:2}
  \end{figure}

\end{example}

Example~\ref{ex:Parken} illustrated how optimal control can be used to simulate a driver in an 
automatic parking maneuver. Now we like to investigate the influence of parameters. 
A parametric sensitivity analysis as in \cite{Mau01} for the nominal optimal control problem OCP($p^*$) 
with some nominal parameter $p^*$ allows to approximate the optimal solution $(z(p),u(p))$ of the 
perturbed problem OCP($p$) with $p$ close to $p^*$ locally by a first order Taylor approximation
\begin{equation} \label{EQ:Sens:1}
  z(p) \approx z(p^*) + \frac{\partial  z}{\partial p}(p^*) (p-p^*),\quad u(p) \approx u(p^*) + \frac{\partial u}{\partial p}(p^*) (p-p^*).
\end{equation}
This approximation holds for $p$ sufficiently close to $p^*$ under suitable assumptions on the nominal solution, i.e. 
first-order necessary conditions, strict complementarity, second-order sufficient conditions and linear independence constraint qualification have to hold. 
A similar approximation holds for discretized optimal control problems, compare \cite{Bue01a}.
These update formulas can be applied in real-time to update a nominal solution 
in the presence of perturbations in $p^*$. This idea can be exploited in multistep model-predictive control schemes 
as well, compare \cite{Palma2015}. 

In the following example, we demonstrate the parametric sensitivity analysis for a collision avoidance maneuver. 

\begin{example}[Collision avoidance maneuver]
  The task is to avoid a collision with a fixed obstacle that is blocking the right half of a straight road. 
  The width of the road is $8$ meters and we aim to find the minimal distance $d$ such that an avoidance maneuver 
  is possible with moderate steering effort. We introduce two perturbation parameters $p_1$ and $p_2$ with nominal values 
  $p_1^* = p_2^* = 0$ into the problem 
  formulation. The first parameter $p_1$ models perturbations in the initial yaw angle. The second parameter $p_2$ models 
  perturbations in the motion of the obstacle and allows the obstacle to move with a given velocity $v_{obs} = 100$ [km/h] into 
  a given direction with angle $\psi_{obs} = 170$ [$^\circ$]. 

  The car's dynamics are given by (\ref{EQ:1a})-(\ref{EQ:1c}), (\ref{EQ:2a}) and (\ref{EQ:2b}) with $\ell=2.7$ and width $b=2$. 
  The steering angle velocity $w$ and the acceleration $a$ are restricted by the control constraints
  \begin{equation}\label{EQ:Avoidance:1}
    w(t) \in [-0.5,0.5],\qquad a(t)\in [-10.0,0.5].
  \end{equation} 
  The steering angle $\delta$ is bounded by the state constraints 
  \begin{equation}\label{EQ:Avoidance:2}
    -\frac{\pi}{6} \leq \delta(t) \leq \frac{\pi}{6}.
  \end{equation}
  The initial state of the car is given by 
  \begin{equation}\label{EQ:Avoidance:3}
    x(0) = 0,\quad y(0) = 1.75,\quad \psi(0) = p_1,\quad \delta(0) = 0, \quad v(0) = 27.78,
  \end{equation}
  where $p_1$ is a perturbation parameter with nominal value $p_1^* = 0$. The terminal state is given by 
  \begin{equation}\label{EQ:Avoidance:4}
    x(t_f) = x_{obs}(t_f) + 3,\quad \psi(t_f) = \delta(t_f) = 0
  \end{equation}
  with 
  \begin{displaymath}
    x_{obs}(t) = d + t p_2 v_{obs} \cos \psi_{obs},\qquad y_{obs}(t) = 3.5 +   t p_2 v_{obs} \sin \psi_{obs}
  \end{displaymath}
  and perturbation parameter $p_2$ with nominal value $p_2^* = 0$. Note that the obstacle is not 
  moving for the nominal parameter $p_2^*=0$. 

  Moreover, the obstacle is defined by the state constraints 
  \begin{equation}\label{EQ:Avoidance:5}
    s(x(t), x_{obs}(t), y_{obs}(t)) + \frac{b}{2} \leq y(t) \leq 8 - \frac{b}{2},
  \end{equation}
  where $s$ is the piecewise defined and continuously differentiable function
  \begin{displaymath}
    s(x,d,h) := \left\{\begin{array}{ll}
     0,  & \mbox{if } x < d,\\
     4 h (x - d)^3, & \mbox{if } d \leq x < d + 0.5, \\
     4 h (x - (d+1))^3 + h, & \mbox{if } d+0.5 \leq x < d+1, \\
     h, & \mbox{if } x \geq d+1.
    \end{array}\right. 
  \end{displaymath}
  
  Summarizing, the optimal control problem aims at minimizing a linear combination of initial distance to the obstacle $d$ and steering effort, that is 
  \begin{displaymath}
    d  + 18 \int_0^{t_f} w(t)^2 dt, 
  \end{displaymath}
  subject to the constraints (\ref{EQ:1a})-(\ref{EQ:1c}), (\ref{EQ:2a}), (\ref{EQ:2b}), (\ref{EQ:Avoidance:1})-(\ref{EQ:Avoidance:5}).   
  
  Figure~\ref{Fig:Ausweichen:1} shows the output of the direct shooting method OCPID-DAE1, see \cite{OCPIDDAE1}, with $N=51$ grid 
  points for the nominal optimal control problem with parameters $p_1^* = p_2^* = 0$. 
  The final time is $t_f\approx 1.00541$, the optimal distance to the obstacle amounts to $d\approx 19.62075$ and the acceleration is active at 
  the lower bound with $a(t) = -10$ for all $t \in [0,t_f]$. 

  \begin{figure}
    \begin{center}
      \includegraphics[height=0.8\textwidth,angle=-90]{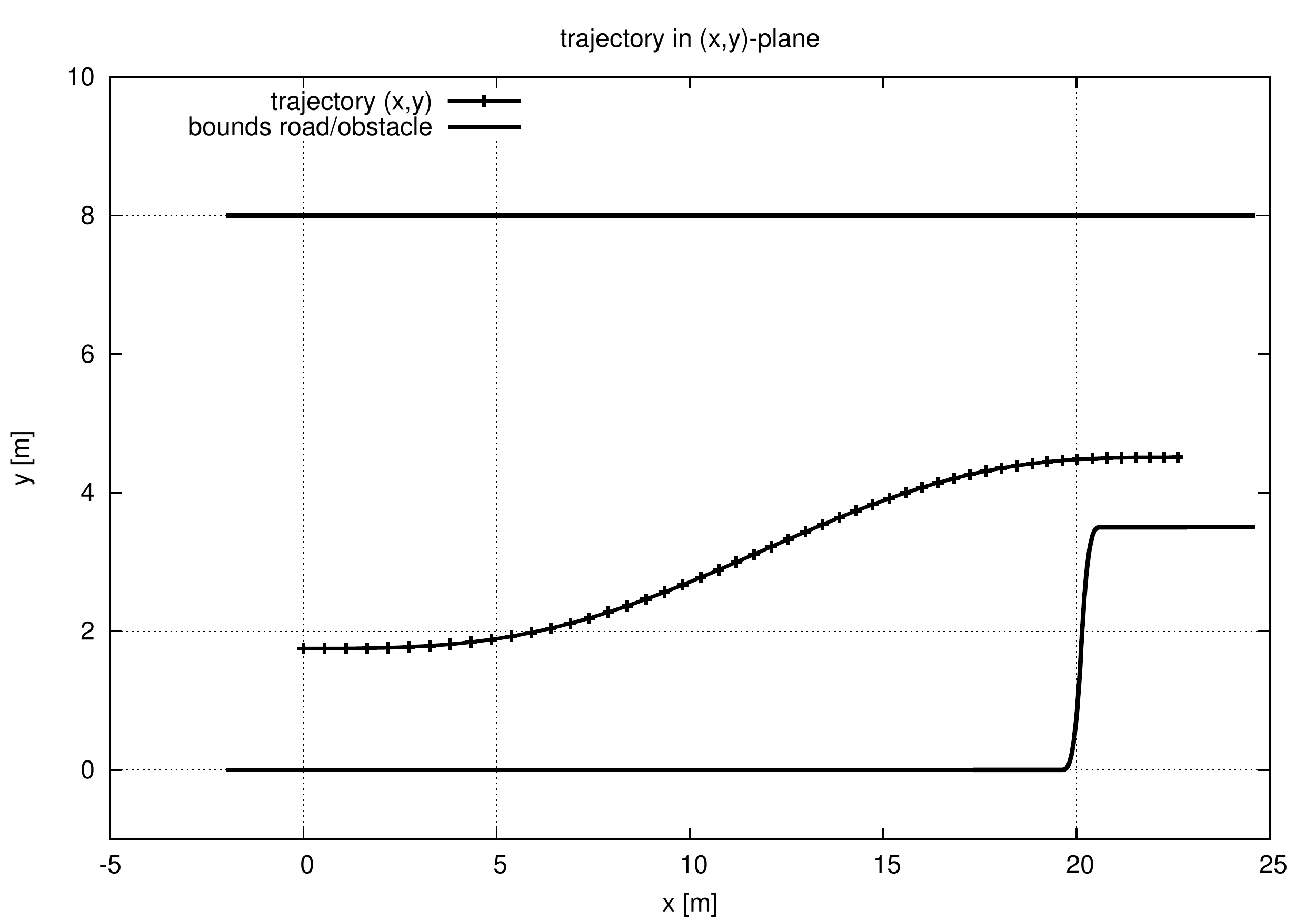}    

      \includegraphics[height=0.3\textwidth,angle=-90]{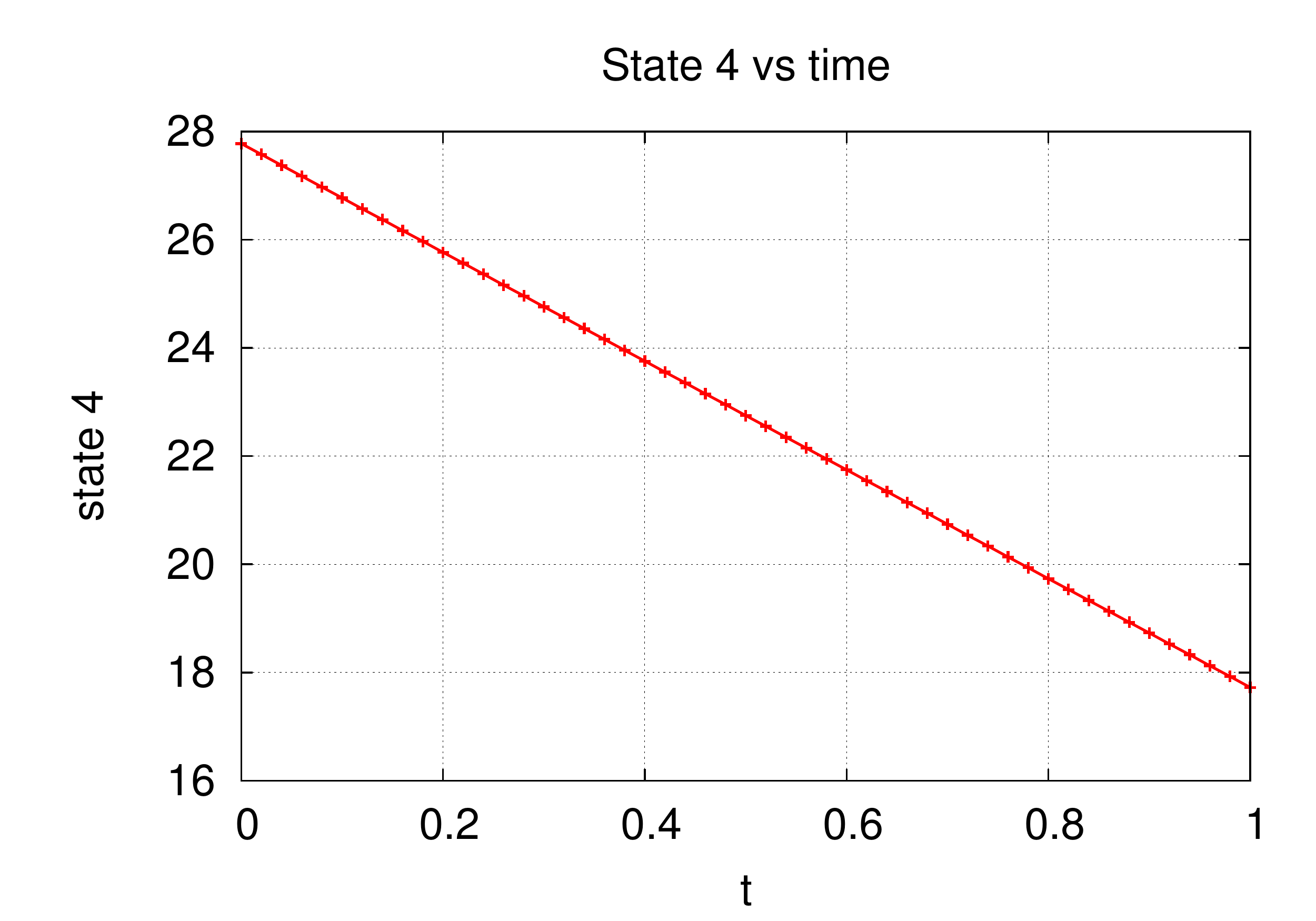}      
      \includegraphics[height=0.3\textwidth,angle=-90]{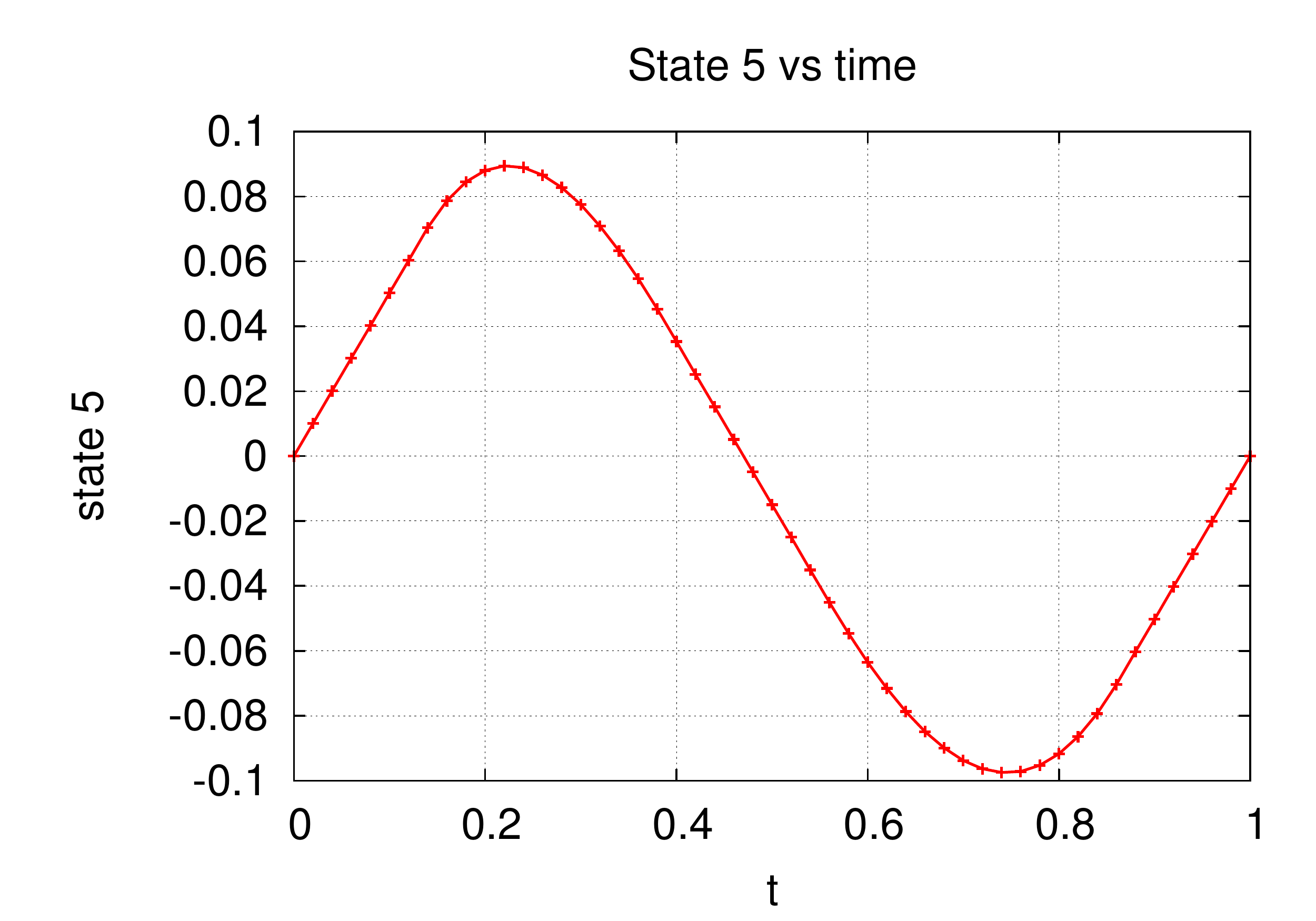}
      \includegraphics[height=0.3\textwidth,angle=-90]{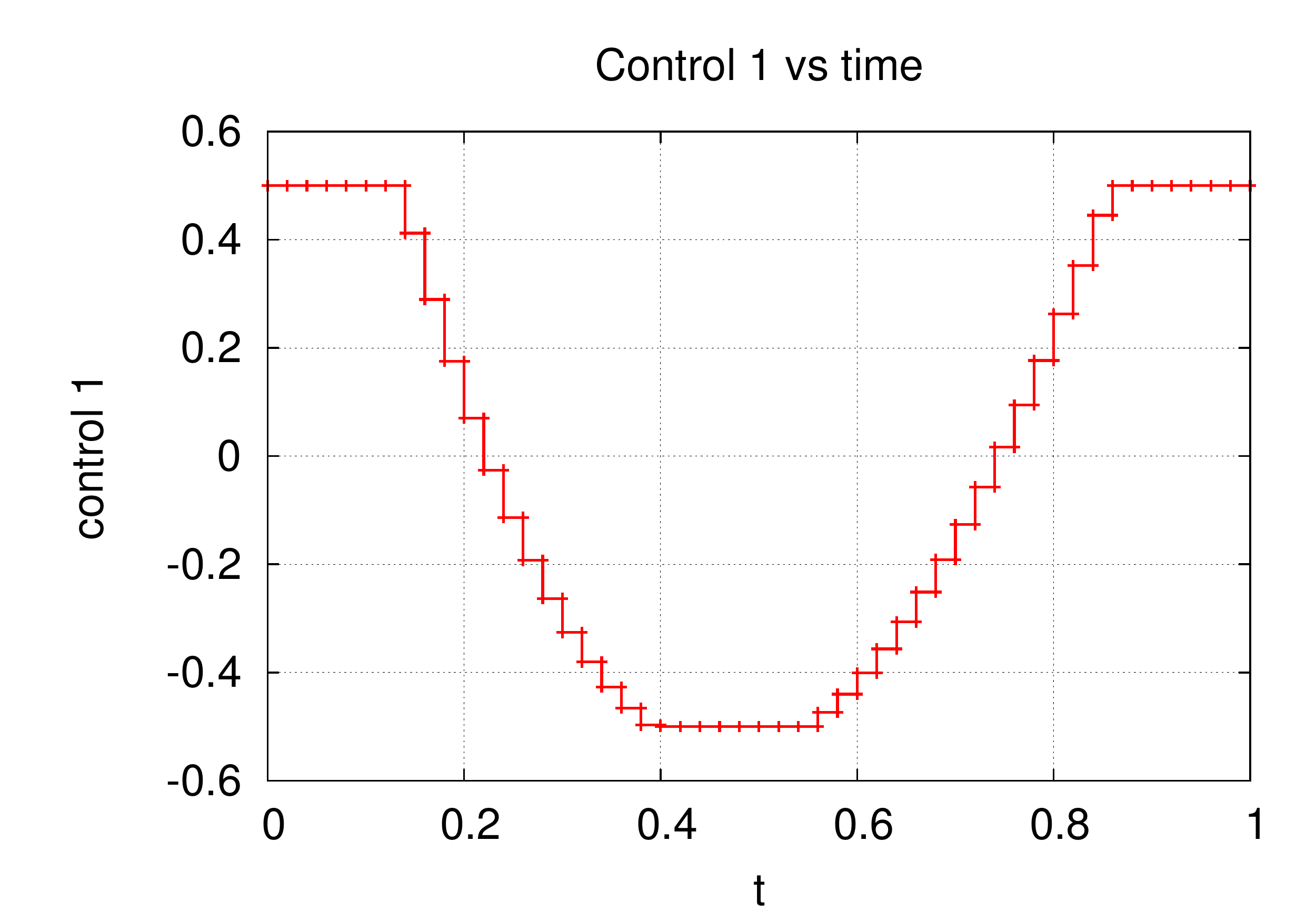}
    \end{center}

    \caption{Trajectory $(x,y)$ (state, top), velocity $v$ (state, bottom left), steering angle $\delta$ (state, bottom middle) 
      and steering angle velocity $w$ (control, bottom right) for optimal avoidance maneuver.}
    \label{Fig:Ausweichen:1}
  \end{figure}
  
  Figure~\ref{Fig:Ausweichen:2} shows the sensitivities $\partial w/\partial p_1$ and $\partial w/\partial p_2$ of the nominal steering angle 
  velocity $w(t)$ with respect to the parameters $p_1$ and $p_2$ for $t\in [0,t_f]$ . 
  
  \begin{figure}
    \begin{center}
      \includegraphics[height=0.8\textwidth,angle=-90]{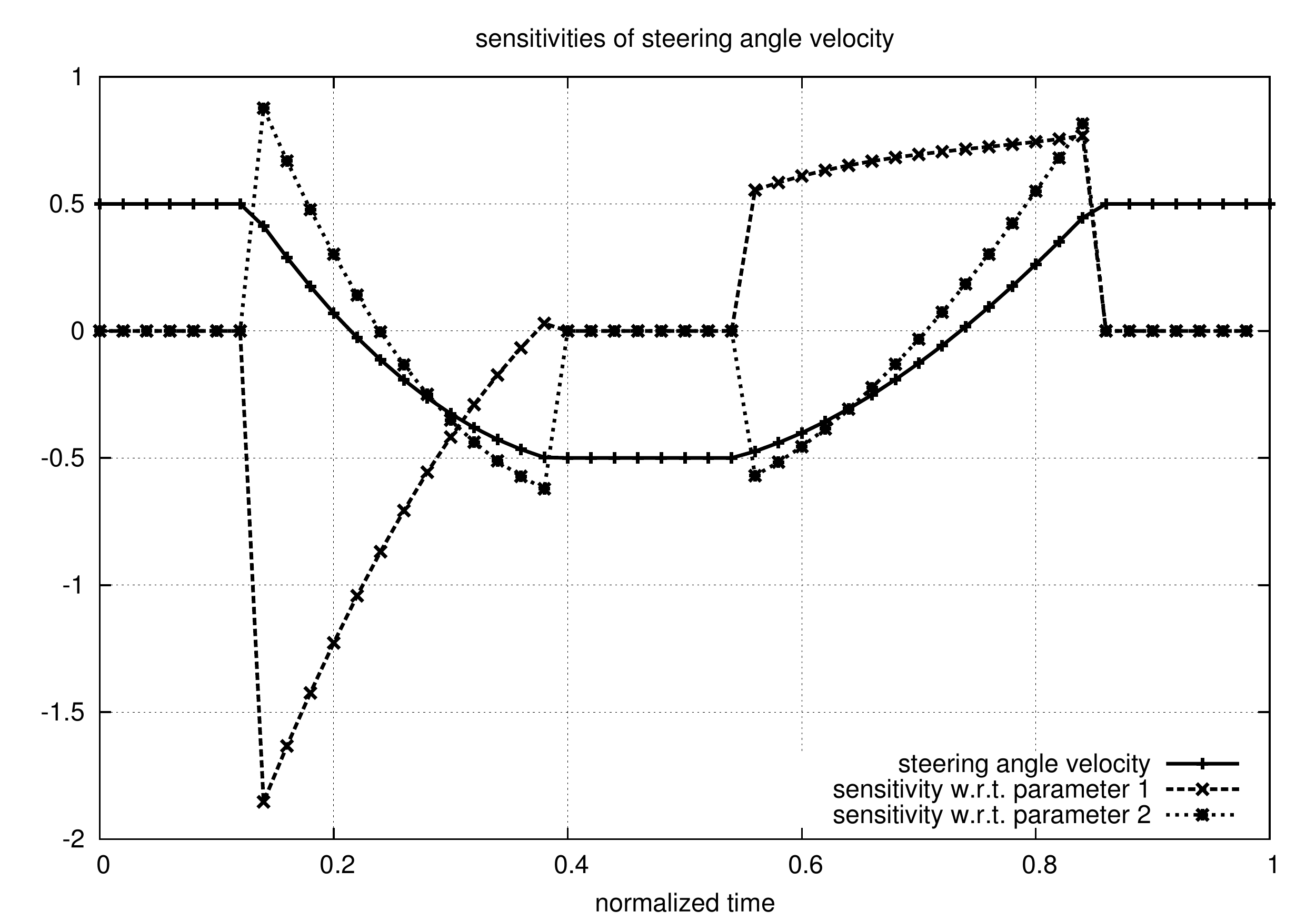}
    \end{center}

    \caption{Sensitivities of the steering angle velocity w.r.t. to $p_1$ (perturbation in initial yaw angle) and 
      $p_2$ (perturbation of obstacle).}
    \label{Fig:Ausweichen:2}
  \end{figure}

  The sensitivities of the nominal final time $t_f$ with respect to $p_1$ and $p_2$ compute to 
  \begin{displaymath}
    \frac{\partial t_f}{\partial p_1} \approx -1.66018,\qquad \frac{\partial t_f}{\partial p_2} \approx 0.50118.
  \end{displaymath}
  The sensitivities of the nominal distance $d$ with respect to $p_1$ and $p_2$ compute to 
  \begin{displaymath}
    \frac{\partial d}{\partial p_1} \approx -28.95949,\qquad \frac{\partial d}{\partial p_2} \approx 35.66225. 
  \end{displaymath}

  The sensitivities allow to predict the optimal solution under (small) perturbations using the Taylor approximation in 
  (\ref{EQ:Sens:1}). Figures~\ref{Fig:Ausweichen:3}-\ref{Fig:Ausweichen:4} show the results of such a prediction for 
  perturbations in the range of $p_1 \in [-0.1,0.1]$ and $p_2\in [-0.1,0.1]$. 

  \begin{figure}
    \begin{center}
      \includegraphics[height=0.8\textwidth,angle=-90]{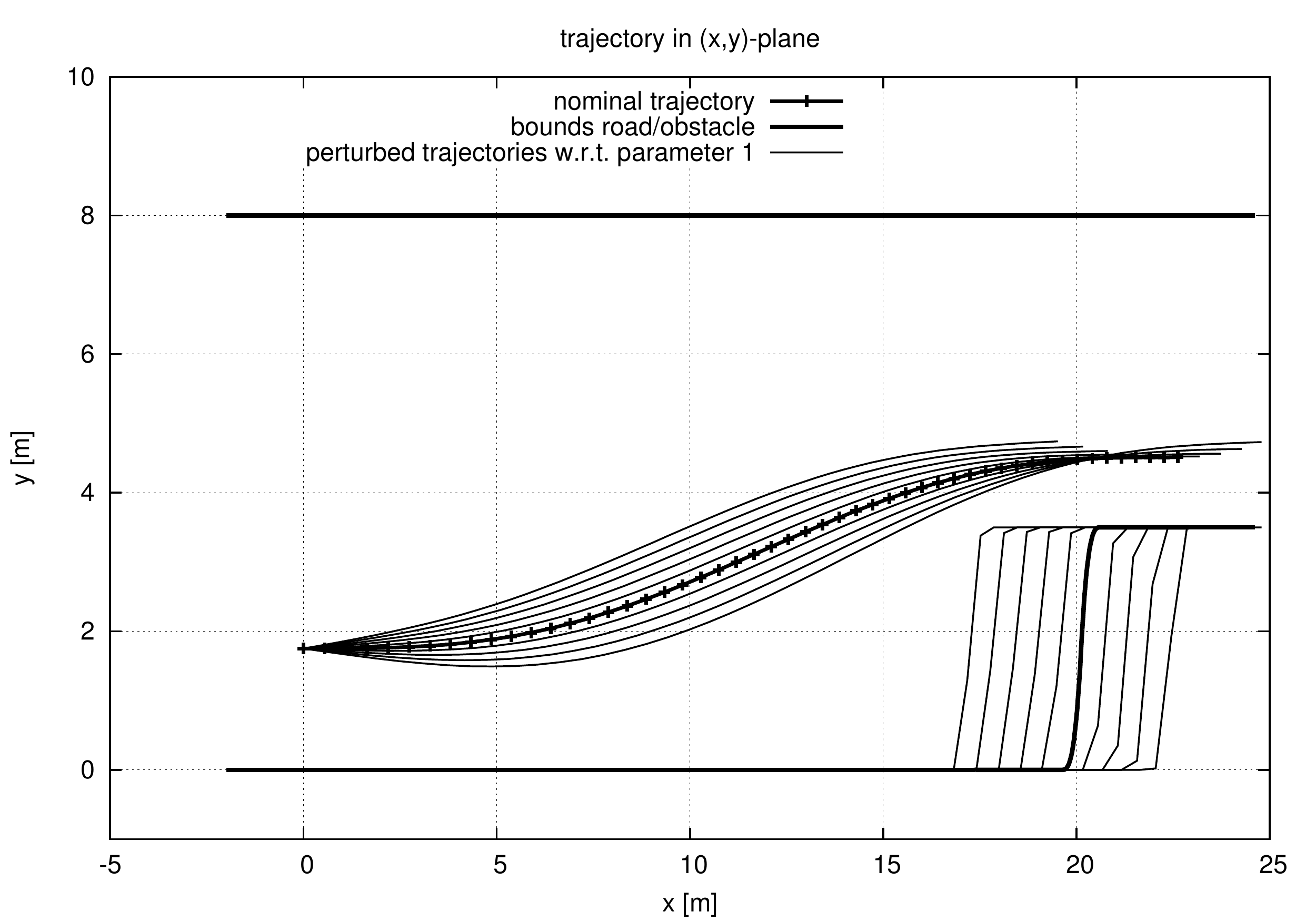} 

      \includegraphics[height=0.4\textwidth,angle=-90]{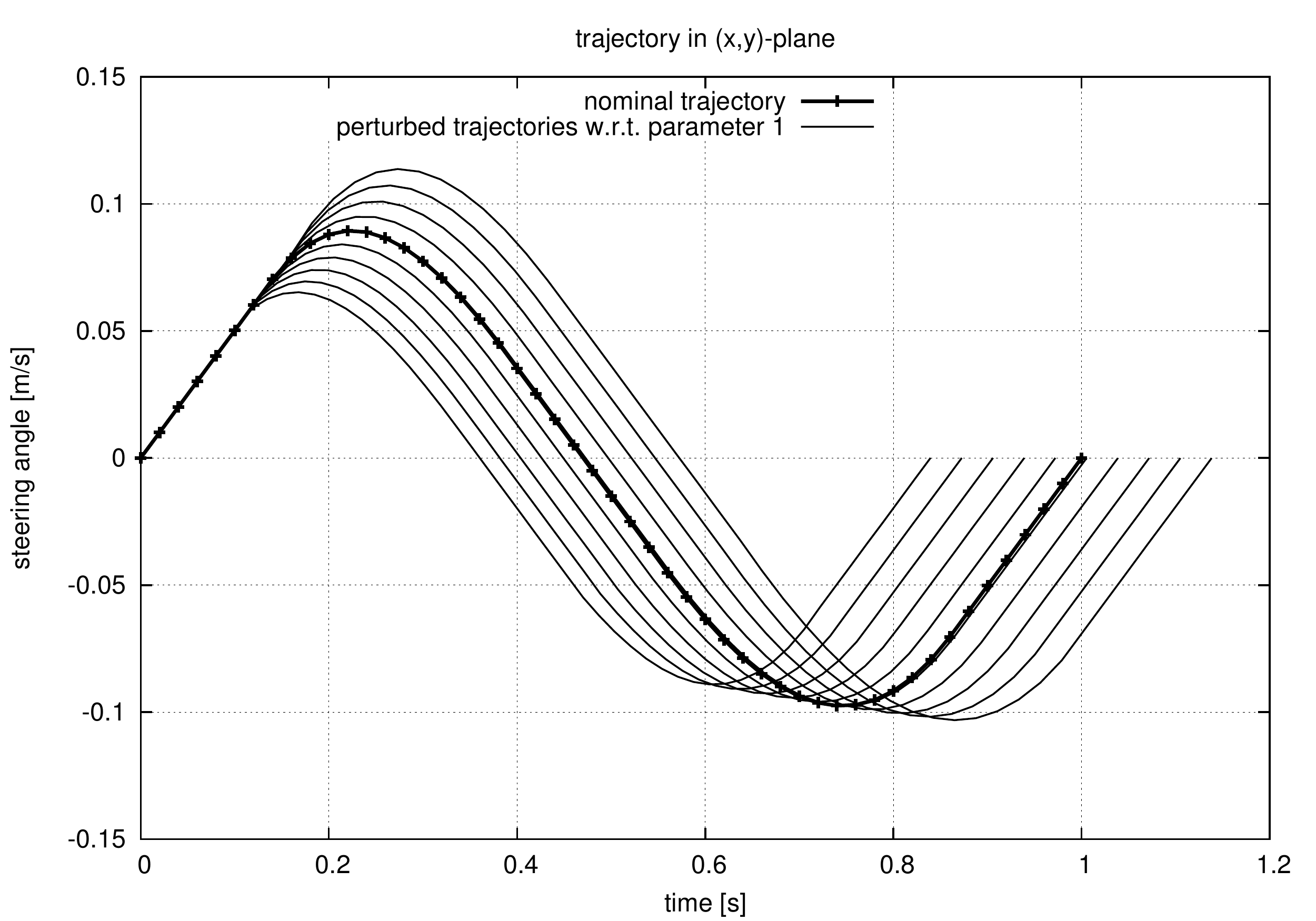} 
      \includegraphics[height=0.4\textwidth,angle=-90]{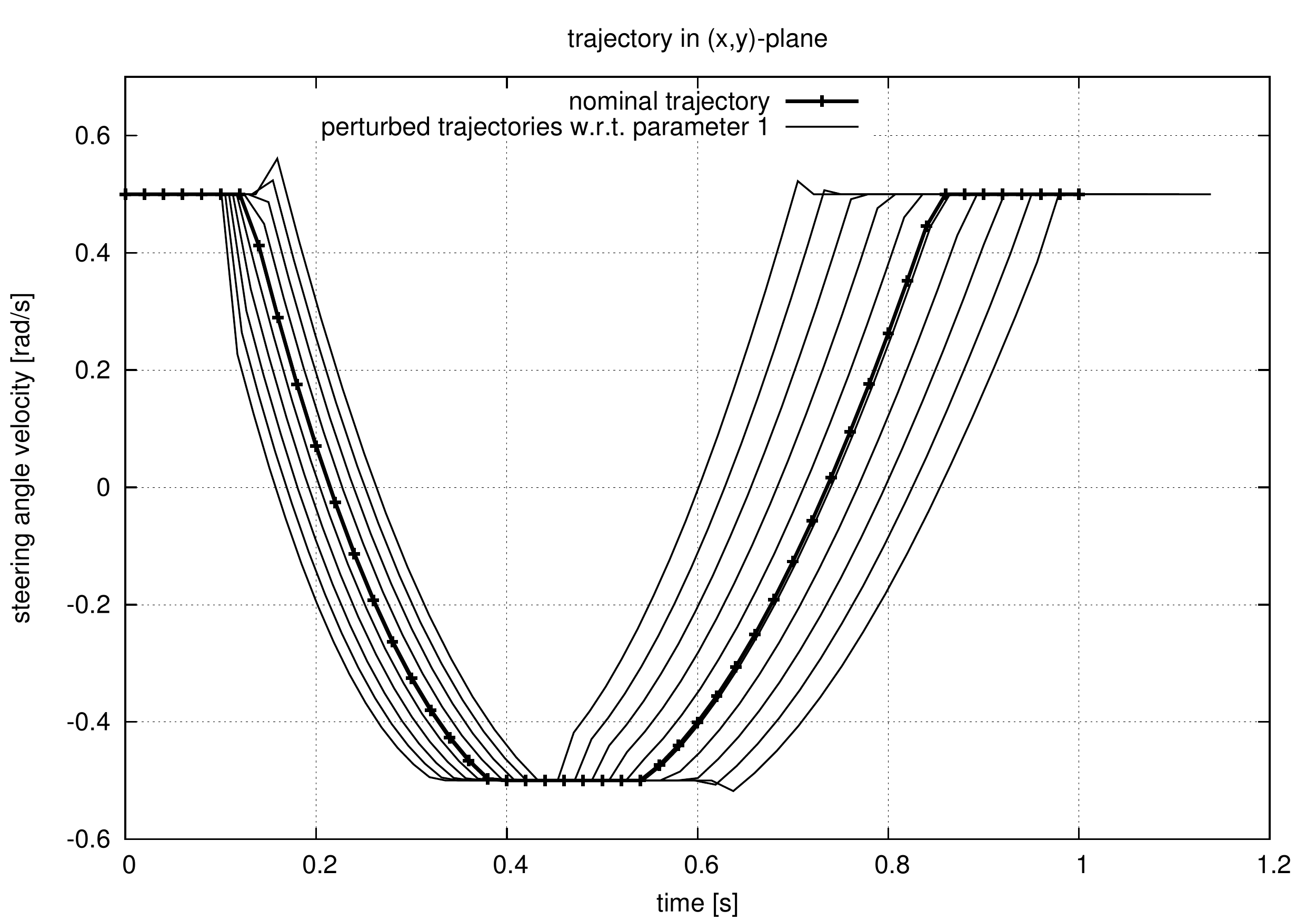} 
    \end{center}

    \caption{Perturbed trajectories (top), steering angle (bottom, left) and steering angle velocity (bottom, right) 
      with sensitivity updates for perturbations $p_1\in [-0.1,0.1]$. Please note that not all perturbations are feasible.}
    \label{Fig:Ausweichen:3}
  \end{figure}
  
  \begin{figure}
    \begin{center}
      \includegraphics[height=0.8\textwidth,angle=-90]{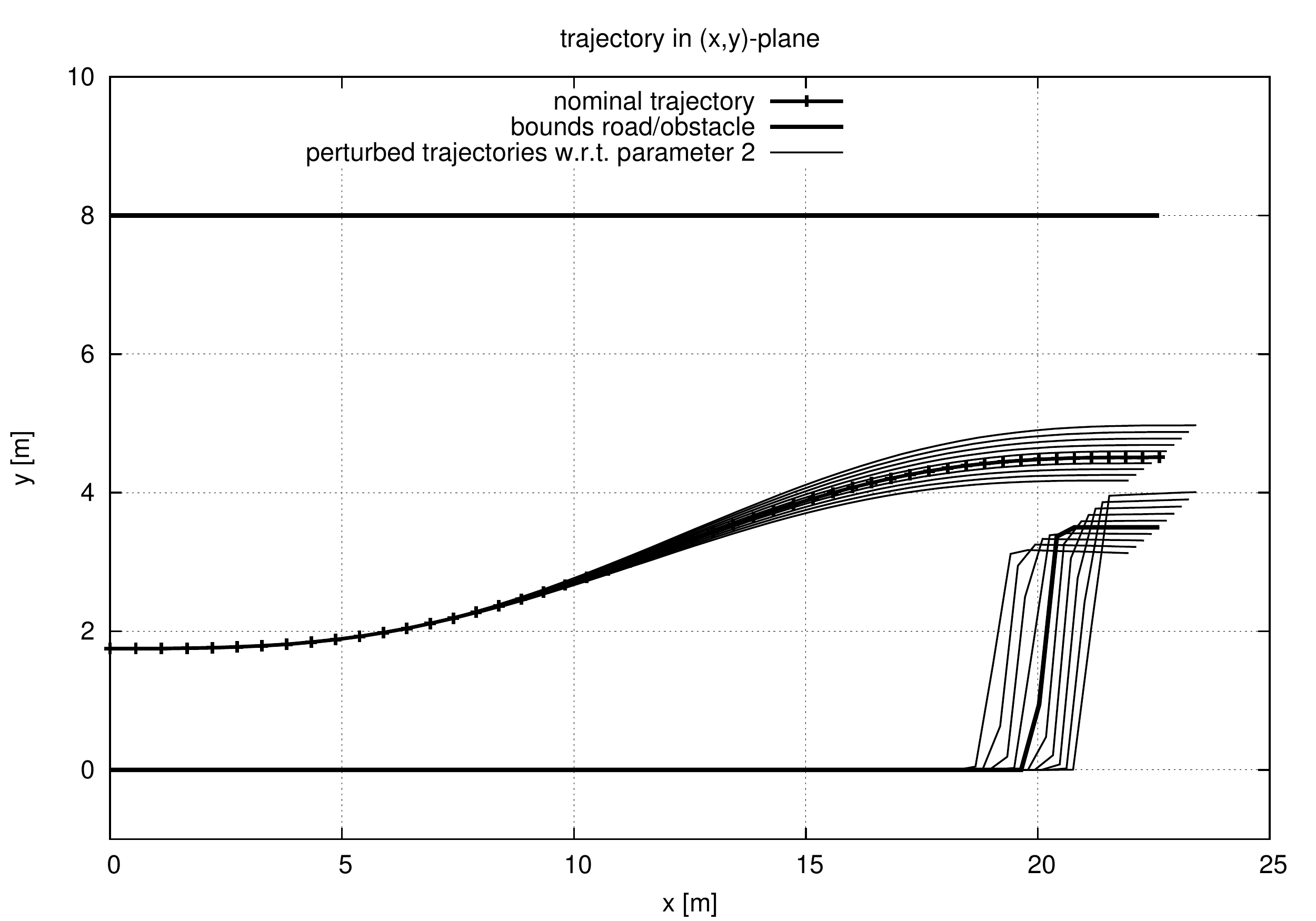} 
      
      \includegraphics[height=0.4\textwidth,angle=-90]{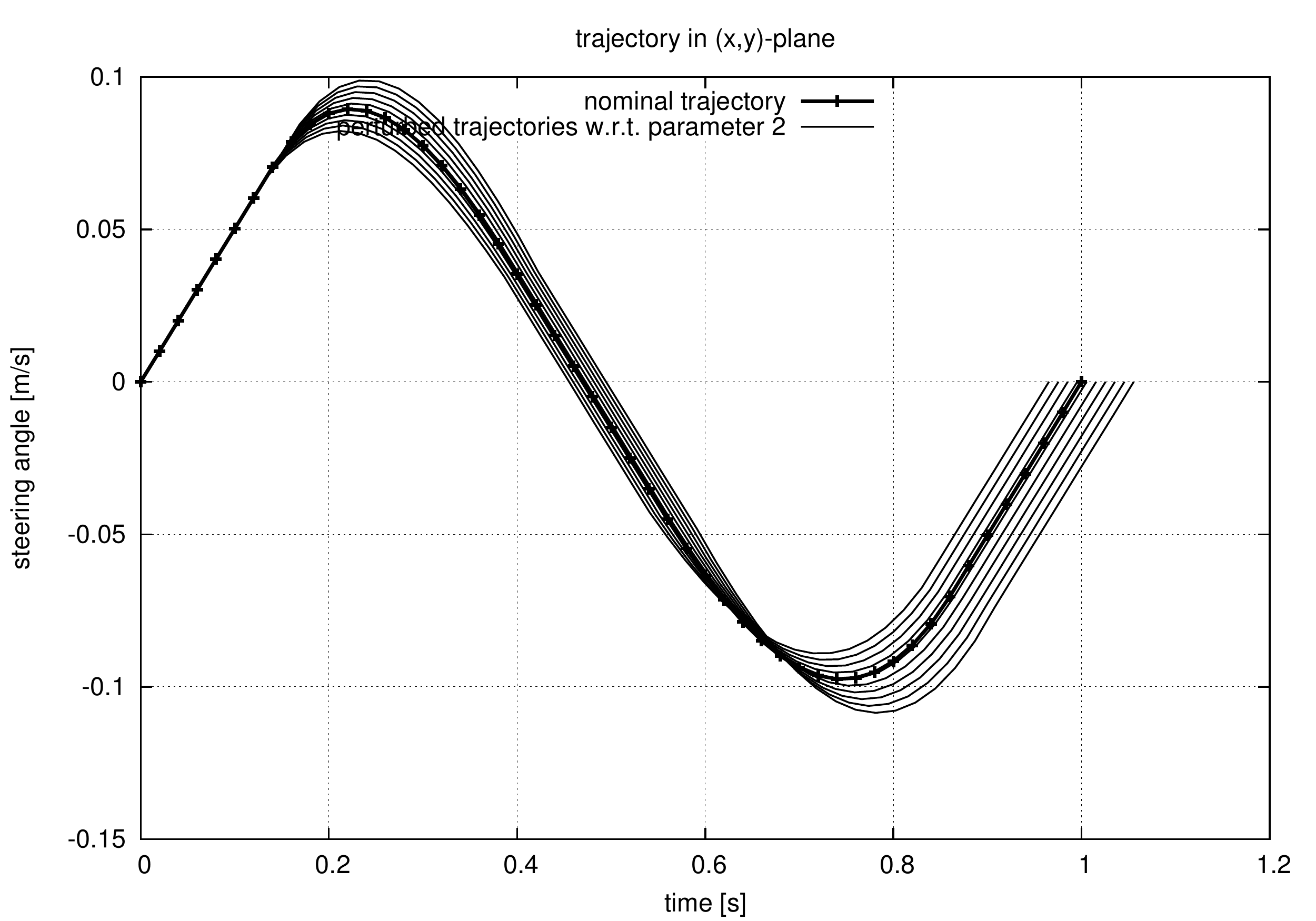} 
      \includegraphics[height=0.4\textwidth,angle=-90]{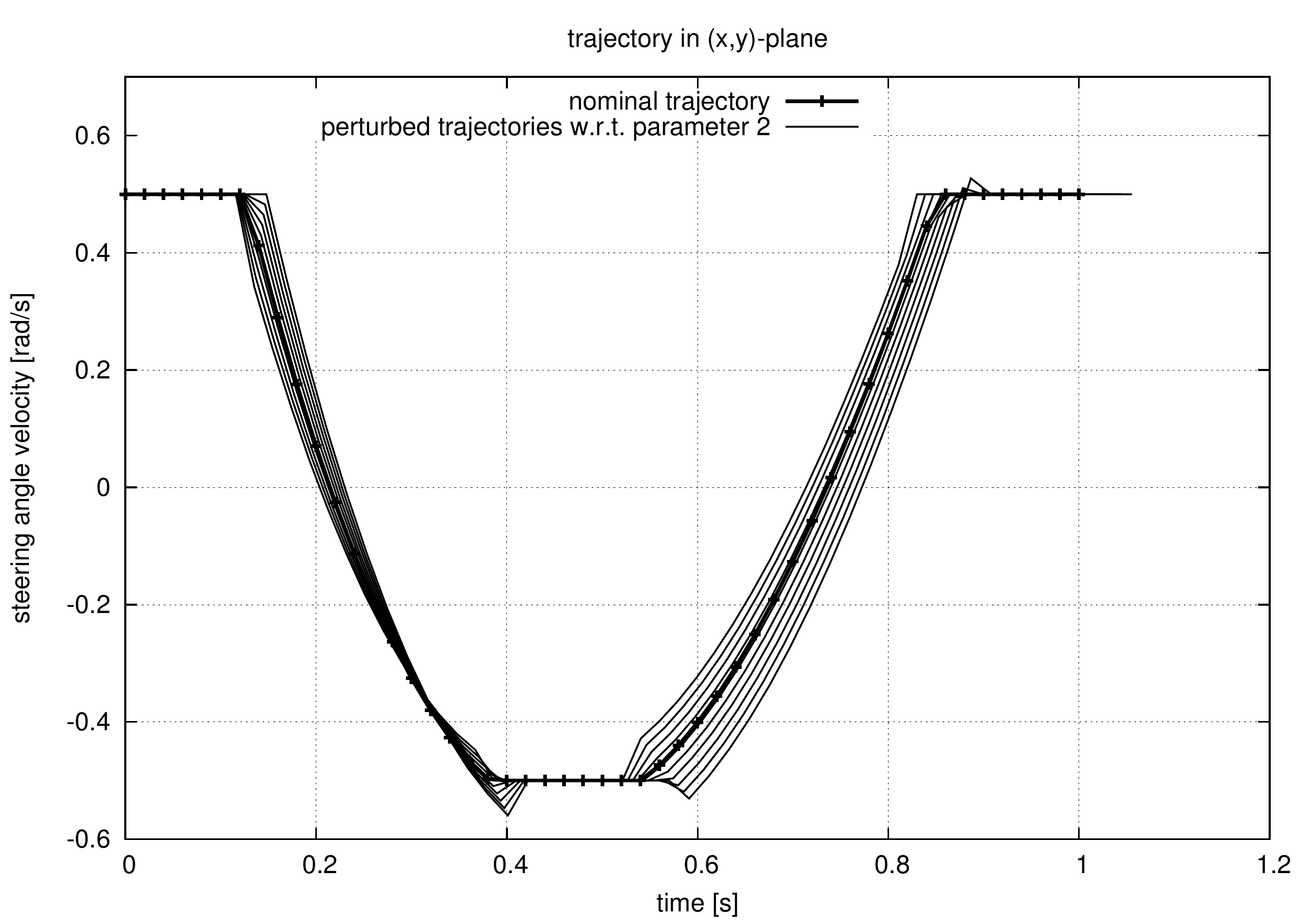} 
    \end{center}
    
    \caption{Perturbed trajectories (top), steering angle (bottom, left) and steering angle velocity (bottom, right) 
      with sensitivity updates for perturbations $p_2\in [-0.1,0.1]$. Please note that not all perturbations are feasible.}
    \label{Fig:Ausweichen:4}
  \end{figure}

\end{example}

Naturally, these examples can only provide a small idea of how optimal control techniques can be 
used to control vehicles and many extensions and applications to more complicated scenarios exist. 
The problem of controlling an autonomous vehicle by means of optimal control in real-time was 
addressed in \cite{Schmidt2014b}. Optimization based obstacle avoidance techniques can be found 
in \cite{Schmidt2014a,Schmidt2014c}.

A different approach that exploits ideas from reachability analysis 
was used in \cite{WooEsf:2012:IFA_4040} to design a controller for a scale car that drives 
autonomously on a given track. Reachable sets turn out to be a powerful tool to detect and avoid 
collisions and to investigate the influence on perturbations on the future dynamic behavior. A comprehensive 
overview can be found in \cite{Kurzhanski2014}. \cite{Althoff2010a} uses reachability analysis with 
zonotopes and linearized dynamics for collision detection. Reachable set approximations 
through zonotopes have been obtained in \cite{Althoff2011f}. Verification approaches for collision avoidance 
systems using reachable sets are investigated in \cite{Nielsson2014,Xausa2015}.

Virtual drives with gear shifts leading to mixed-integer optimal control problems 
have been considered in \cite{Ger05a,Ger06b,Kirches2010}.

\section{Distributed Hierarchical Model-Predictive Control for Communicating Vehicles}
\label{sec:control}

While the trajectory generation for a single vehicle was in the focus of Section~\ref{sec:trajectory}, 
we are now discussing control strategies for several autonomous vehicles that interact with each 
other. To this end we assume that we have $N\in\N$ vehicles that can communicate through suitable communication 
channels and exchange information on positions, velocities, and predicted future behavior. The aim is to efficiently control 
the vehicles in a self-organized and autonomous way without prescribing a route. We suggest to use a distributed 
model-predictive control (MPC) strategy and couple it with a priority list or hierarchy, compare 
\cite{Kianfar2012,Pannek2013,Gross2014,Findeisen2014}. The priority list will rank the vehicles in an adaptive way 
depending on the current driving situation and give highly ranked vehicles priority while driving. 
Vehicles with low priority have to obey the motion of vehicles with higher priority. 

\begin{figure}[h]
\setlength{\unitlength}{1pt} \centering
  \begin{picture}(250,130)(0,-20)
    \put(-40,-10){\includegraphics[width=0.5\textwidth]{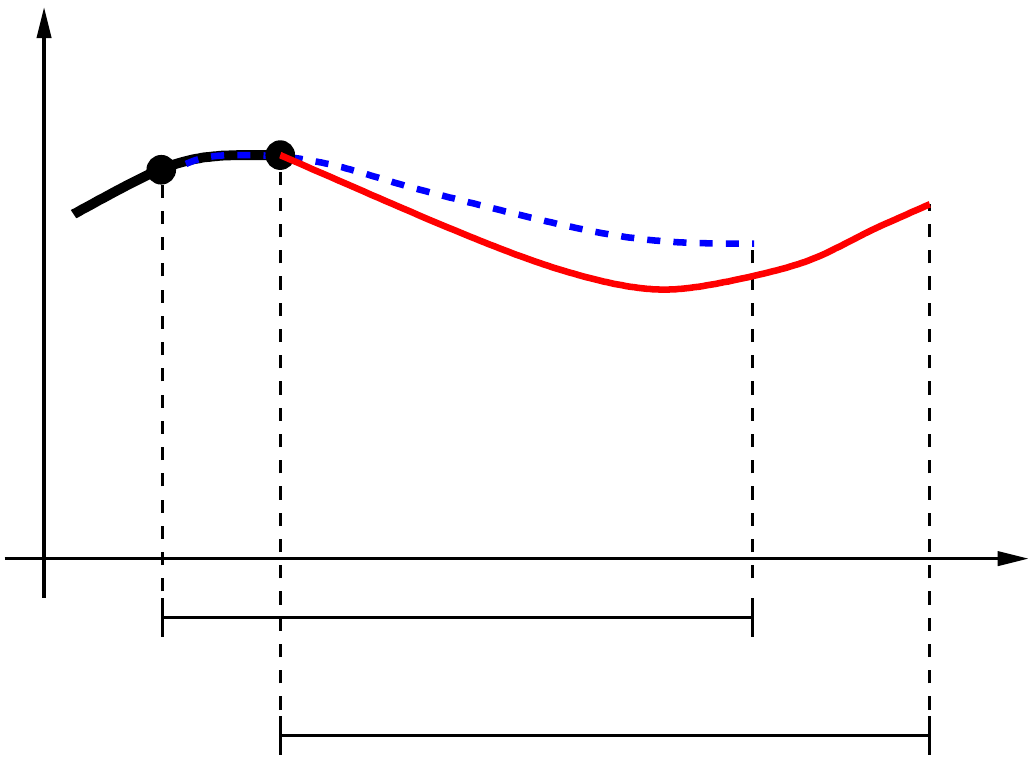}}
    \put(-15,1){$t_i$}
    \put(75,1){$t_i+T$}
    \put(0,-20){$t_{i+1}=t_i+\tau$}
    \put(100,-20){$t_{i+1}+T$}
    \put(40,85){$\hat z\big|_{[t_i,t_i+T]}$}
    \put(30,58){$\hat z\big|_{[t_{i+1},t_{i+1}+T]}$}
    \put(-10,95){$z(t)$}
    \put(115,115){\includegraphics[height=0.55\textwidth,angle=-90]{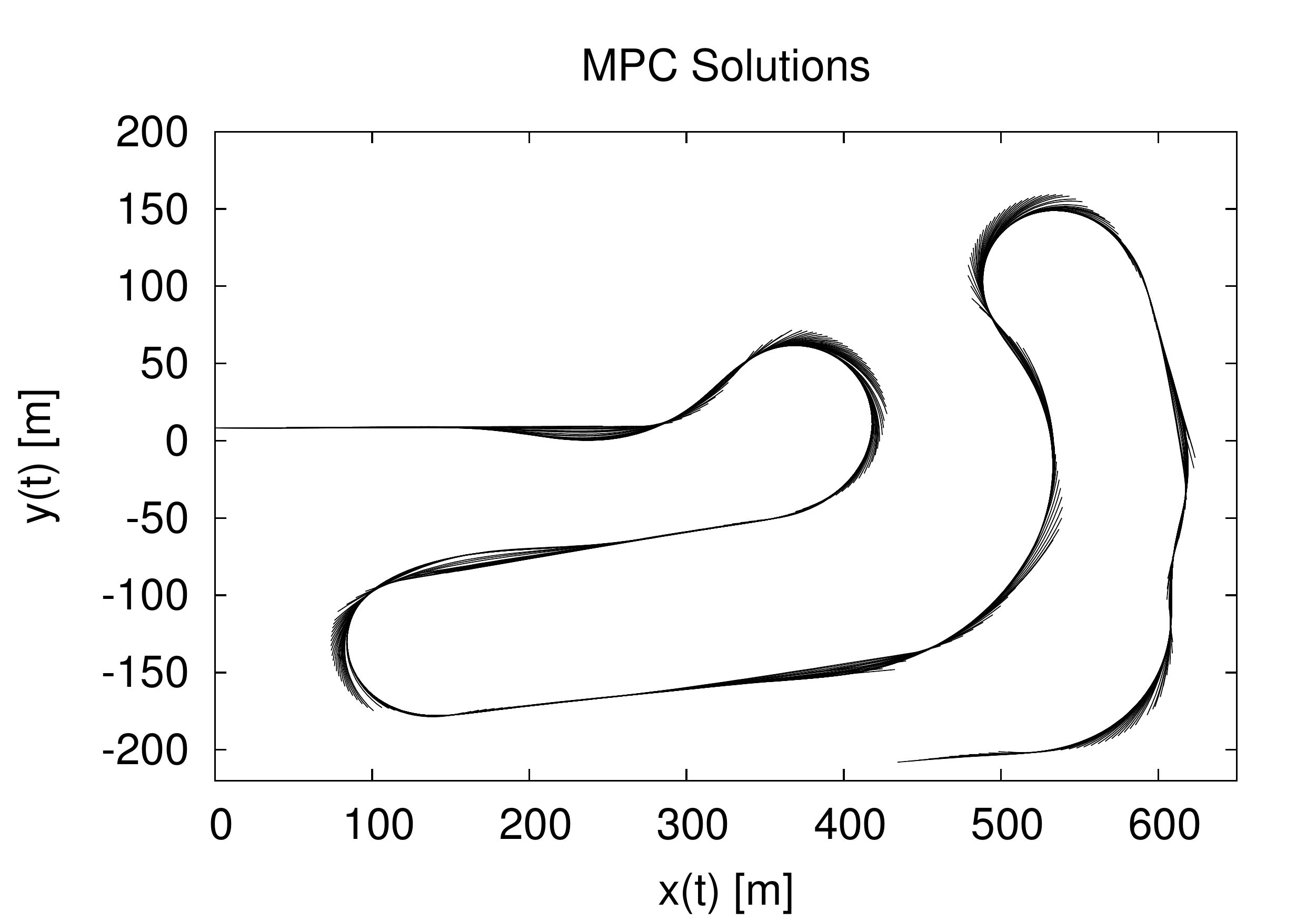}}
    \end{picture}
  \caption{Model-predictive control scheme: Repeated optimization on a prediction horizon of length $T$ and 
    control acting on a control interval of length $\tau$ (left). Application to a test-drive (right): Local solutions
    of the MPC scheme, compare \cite{Ger09}.}\label{mg:FIG:mpc_concept} 
\end{figure}

Model-predictive control is a well established feedback control paradigm, compare \cite{Gru11} for a detailed 
exposition and discussion of stability and robustness properties. The working principle of model predictive 
control is based on a repeated solution of (discretized) optimal control problems on a moving time horizon, 
see Figure~\ref{mg:FIG:mpc_concept}. The model-predictive control scheme was used in \cite{Ger09} to simulate the 
drive on long tracks, see the picture on the right in Figure~\ref{mg:FIG:mpc_concept} for an example. 

The MPC algorithm depends on a local time horizon of length $T>0$, sampling times $t_i$, $i\in\N_0$, and 
a shifting parameter $\tau > 0$. On each local time horizon $[t_i,t_i+T]$ 
an optimal control problem has to be solved for a given initial state, which may result from 
a measurement and represents the current state of the vehicle. 
Then the computed optimal control is implemented on the interval $[t_i,t_i+\tau]$ 
and the state is measured again at the sampling time $t_{i+1} = t_i+\tau$. 
After that, the process is repeated on the shifted time interval $[t_{i+1},t_{i+1}+T]$. This control paradigm provides a 
feedback control, since it reacts on the actual state at the sampling times $t_i$, $i\in\N_0$. 
Moreover, the control paradigm is very 
flexible since control and state constraints and individual objectives can be considered in the local optimal control 
problems, compare Section~\ref{sec:OCP}.

The basic model-predictive control scheme has to be enhanced towards distributed vehicle systems. Herein, 
the computations take place on the individual vehicles and the relevant information is exchanged. 
In addition, a priority list has to be included.

\subsection{Priority List}
\label{sec:Priority}

The priority list consists of a set of predefined rules to rank the 
vehicles. Vehicles with a lower priority have to take into account in their motion planning algorithm
the motion of the vehicles with higher priority, while vehicles with the same priority 
can move independently, see  \cite{Pannek2013}.

To simplify the computations we assume that each vehicle only considers its neighboring vehicles as potential obstacles 
in the state constraints, i.e. vehicles which are inside a certain communication radius or distance. 
Vehicles outside this neighborhood are ignored in the optimization process. For each of the $N\in\N$ vehicles we introduce the set
\[I_{NH}(i):=\left\lbrace j \in\left\lbrace 1,\ldots,N \right\rbrace\setminus\left\lbrace i\right\rbrace \mid j \text{ is a neighbor of } i \right\rbrace, \]
which contains the indexes of the neighbors of the $i$-th vehicle. 

In order to avoid a conflict between vehicles, which could lead to a collision, we introduce a set of rules which assigns a distinct hierarchy level to each vehicle, i.e. every 
vehicle inside the $i$-th vehicle's neighborhood has either a higher or lower hierarchy level than the $i$-th vehicle. The vehicle with the lower hierarchy level has to consider the safety boundary of the vehicles with a higher hierarchy level in terms of state constraints, while a vehicle with a higher priority is allowed to ignore the safety boundaries of neighboring vehicles with a lower priority. The set 
\[I_{PR}(i):=\left\lbrace j \in I_{NH}(i)\mid j \text{ has a higher priority than } i \right\rbrace \]
contains all indexes of the neighboring vehicles of the $i$-th vehicle, which have a higher hierarchy level than the $i$-th vehicle. If a vehicle has the highest possible priority this set is empty, if it has the lowest hierarchy level it contains all neighboring vehicles. This approach is also able to handle vehicles, which are not part of the communication network, e.g. an ambulance or non autonomous vehicles, by giving those vehicles the highest priority. The rules may also have a certain priority, e.g. traffic rules may have a higher priority than mathematically motivated rules. 

By this approach the computational effort is being reduced, because vehicles with no neighbors or vehicles with the highest hierarchy level are allowed to ignore other vehicles and are able to drive in an optimal way with fewer state constraints. It also reduces the potential for conflicts between vehicles, because of the distinct hierarchy. Depending on the scenario, the rules which assign the hierarchy level might change, e.g. if two vehicles meet at an intersection the vehicle coming from the right would have higher priority, but in a roundabout scenario the vehicle inside the roundabout would have the highest hierarchy level. To identify the priority we first considered traffic rules in our analysis. If those rules fail to assign a distinct priority we then used a mathematically motivated rule which is based on the adjoints computed by solving the optimal control problem. By using the adjoints for the controls we determined which vehicle would have the higher cost if it would deviate from its optimal trajectory. For a fixed set of priority rules we get the following 
distributed hierarchical model-predictive control algorithm:\\
\\
\begin{tabular}{|l|}
\hline 
\\
\Large{\textbf{Distributed Hierarchical MPC Algorithm}} \\ 
\\
\hline 
\\
\textbf{Input:} prediction horizon, control horizon, set of priority rules\\ 
\\
\hline 
\\
\begin{minipage}{1.0\textwidth}

\begin{description}
\item[\textbf{1. }] Determine current states of all vehicles. \\
\item[\textbf{2a.}] Compute in parallel the optimal driving paths of all vehicles with respect to the neighborhood relations and hierarchy levels. \\
\item[\textbf{2b.}] For all vehicles 
\begin{description}
\item[(i)  ] reset all previous neighborhood relations.
\item[(ii) ] screen for neighboring vehicles.
\item[(iii)] submit current states and optimal driving paths to all neighbors.
\end{description}
\item[\textbf{2c.}] For all vehicles 
\begin{description}
\item[(i)] reset all previous priority relations.
\item[(ii)] apply the priority rules and assign the appropriate hierarchy levels.
\end{description}
\item[\textbf{3. }] Apply the computed optimal control on the given control horizon and repeat on shifted time horizon.
\end{description}

\end{minipage}

\\
\\
\hline
\end{tabular} 

\bigskip
The computation of the optimal trajectory for a vehicle stops if the vehicle is close to the destination or if a fixed time limit is exceeded. 
The advantage of a model-predictive control approach is that after each iteration it is possible to update the neighborhood 
relations and the hierarchy levels. 

\subsection{Optimal Control Problems on Prediction Horizon}
\label{sec:OCP-MPC}

In Step 2a of the distributed hierarchical MPC algorithm each of the $N$ vehicles has to solve an individual optimal control problem of type 
\eqref{EQ:OCP:OBJ} - \eqref{EQ:OCP:4}, whose details will be defined in this section. We again use the dynamics 
(\ref{EQ:1a})-(\ref{EQ:1c}), (\ref{EQ:2a}) and (\ref{EQ:2b}) with control vectors $u^{[j]}=(w^{[j]},a^{[j]})^\top$ and state vectors 
$z^{[j]}=(x^{[j]},y^{[j]},\psi^{[j]},v^{[j]},\delta^{[j]})^\top$ for the vehicles $j=1,\ldots,N$. 
We assume that the MPC scheme has proceeded until the sampling point $t_i$ with states $\bar z^{[j]} = z^{[j]}](t_i)$, 
$j=1,\ldots,N$. Moreover, we assume that each vehicle $j\in\{1,\ldots,N\}$ has a given target position 
$(x_\star^{[j]},y_\star^{[j]})$ that it aims to reach. Let $I_{PR}(j)(t_i)$, $j\in \{1,\ldots,N\}$, denote the 
priority sets of the vehicles at time $t_i$. Then the j-th vehicle has to solve the following optimal control 
problem on the local time horizon $[t_i,t_i+T]$ in Step 2a:

\medskip
{\em Minimize
  \begin{align}\label{EQ:MPC:OBJ}
J(z^{[j]},u^{[j]}) &= \alpha_1\,\left\Vert\left(
\begin{array}{c}
x^{[j]}(t_i+T) \\ 
y^{[j]}(t_i+T)
\end{array} 
\right) - 
\left(
\begin{array}{c}
x^{[j]}_\star \\ 
y^{[j]}_\star
\end{array}
\right)\right\Vert^2 \nonumber \\ 
 & \quad + \alpha_2\,\int_{t_i}^{t_i+T} a^{[j]}(t)^2dt + \alpha_3\,\int_{t_i}^{t_i+T} w^{[j]}(t)^2dt
\end{align}
\indent subject to the constraints (\ref{EQ:1a})-(\ref{EQ:1c}), (\ref{EQ:2a}), (\ref{EQ:2b}) 
for $z^{[j]}$ and $u^{[j]}$, $z^{[j]}(t_i) =\bar z^{[j]}$  and 
\begin{displaymath}
  v^{[j]}(t) \in [v_{min}^{[j]}, v_{max}^{[j]}], \quad 
  a^{[j]}(t) \in [a_{min}^{[j]},a_{max}^{[j]}], \quad 
  w^{[j]}(t) \in [-w_{max}^{[j]},w_{max}^{[j]}], 
\end{displaymath}
\indent and
\begin{align*}
  \left(\begin{array}{c} x^{[j]}(t) \\ y^{[j]}(t) \end{array}\right) & \in \Omega_{r} \cap \Omega_{c}^{[j]}(t)
\end{align*}
\indent for $t\in [t_i,t_i+T]$. 
}

\medskip 
Herein, the set $\Omega_r \subseteq \R^2$ defines state constraints imposed by the road. The set $\Omega_{c}^{[j]}(t)$ 
defines time dependent collision avoidance constraints imposed by vehicles with higher priority level, i.e. 
\begin{displaymath}
  \Omega_{c}^{[j]}(t) = \bigcap_{k\in I_{PR}(j)(t_i)} \left\{ \left(\begin{array}{c} x \\ y \end{array}\right)\in \R^2 \;\Bigg|\; \left(\begin{array}{c} x - x_c^{[k]}(t) \\ y - y_c^{[k]}(t) \end{array}\right)^\top 
  Q^{[j]}(\psi^{[k]}(t)) \left(\begin{array}{c} x - x_c^{[k]}(t) \\ y - y_c^{[k]}(t) \end{array}\right) \geq 1 \right\}.
\end{displaymath}
Herein, we assumed an ellipsoidal shape of the vehicles $k\in I_{PR}(j)(t_i)$ with half radii $r_x^{[k]}>0$, $r_y^{[k]} > 0$,  matrix
\begin{displaymath}
  Q^{[k]}(\psi) := S(\psi) \left(\begin{array}{cc} \frac{1}{\left(r^{[k]}_x\right)^2} & 0 \\ 0 & \frac{1}{\left(r^{[k]}_y\right)^2} \end{array}\right)
  S(\psi)^\top
\end{displaymath}
with the rotation matrix $S$ from (\ref{EQ:Drehmatrix}) and the vehicle's center 
\begin{displaymath}
  \left(\begin{array}{c} x^{[k]}_c(t) \\ y^{[k]}_c(t) \end{array}\right) = \left(\begin{array}{c} x^{[k]}(t) \\ y^{[k]}(t) \end{array}\right)
  + S(\psi^{[k]}(t)) \left(\begin{array}{c} \ell^{[k]}/2 \\ 0 \end{array}\right),
\end{displaymath}
where $\ell^{[k]}$ is the length of vehicle $k$, compare (\ref{EQ:CarCenter}).

The weights $\alpha_1^{[j]},\alpha_2^{[j]},\alpha_3^{[j]}>0$ can be used to individually weight the terms in the objective function
(\ref{EQ:MPC:OBJ}) in order to model different drivers. The main goal of the vehicle is to reach its destination, i.e. to 
minimize the distance to its final destination 
$(x_\star^{[j]},y_\star^{[j]})$. The optimal control problem also allows to consider criteria such as fuel consumption or 
comfort in the minimization process by choosing moderate weights $\alpha_2^{[j]}$ and $\alpha_3^{[j]}$, which represent the cost for 
accelerating/braking and steering, respectively.

\subsection{Simplifications}

\begin{figure}[h]
  \begin{center}
    \includegraphics[width=0.8\textwidth]{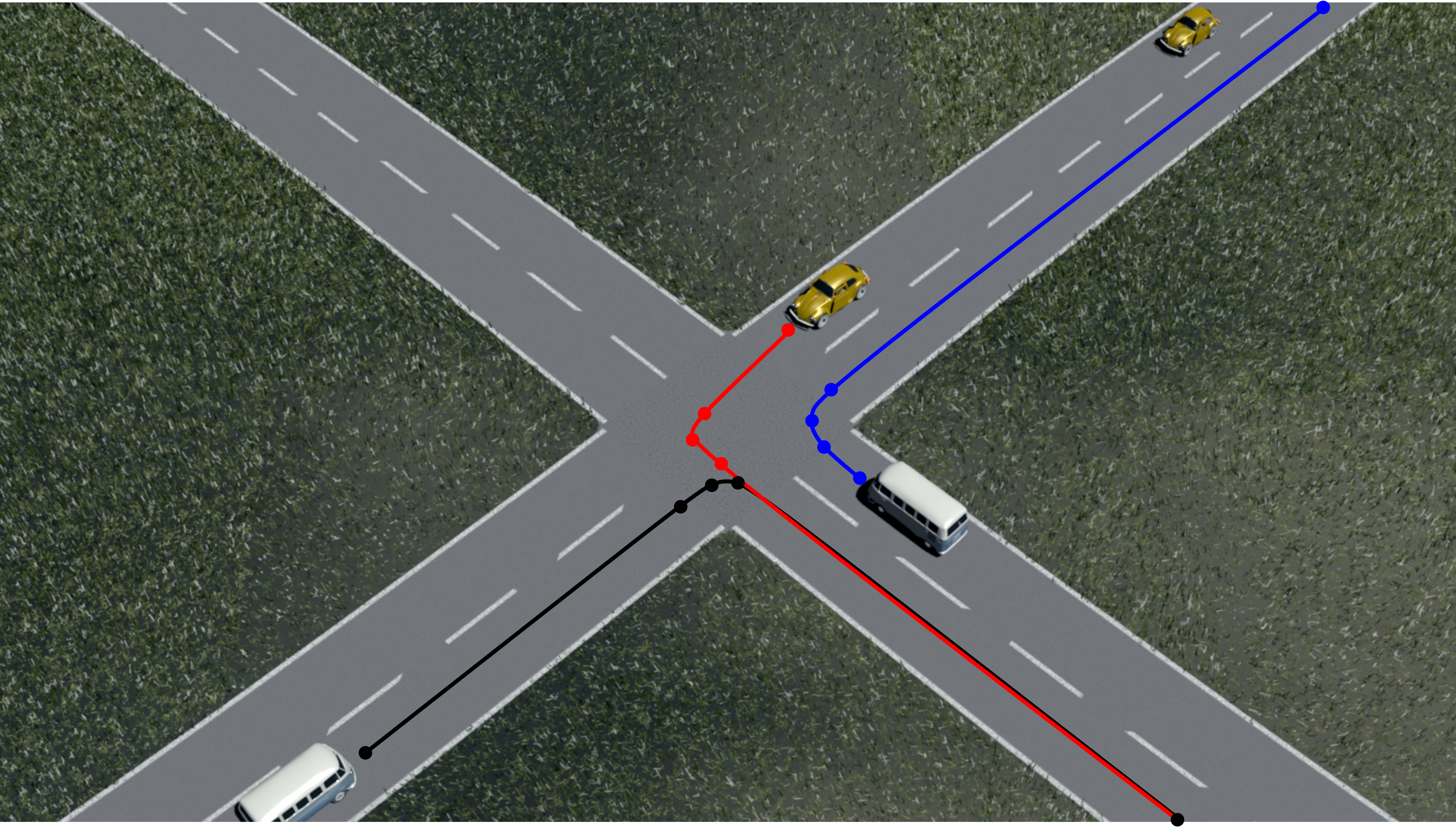}
  \end{center}
  \caption{Distributed hierarchical control of autonomous cars along preassigned paths.}\label{fig:Spline}
\end{figure}

The approach in Section~\ref{sec:OCP-MPC} provides full flexibility to the motion of the vehicles as long as they stay on the 
road and obey collision avoidance constraints. As a result the local optimal control problems are very nonlinear and 
occasionally require a high computational effort, especially if many vehicles interact. One way to decrease the computational 
effort is to simplify the car model so that the cars follow precomputed feasible trajectories coming, e.g., from a 
navigation system, compare \cite{Geyer2014} and Figure~\ref{fig:Spline}. By restricting the vehicle's motion to a predefined 
curve, the degrees of freedom in the local optimal control problems are reduced considerably and real-time computations 
become realistic at the cost of reduced maneuverability.

Let the trajectories of the $N$ vehicles be defined by cubic spline curves 
\[\gamma^{[j]}(s) = \left(\begin{array}{c}
  x^{[j]}(s) \\ 
  y^{[j]}(s)
\end{array}\right),\qquad 0\leq s \leq L^{[j]}, j=1,\ldots,N,\]
which interpolate given way-points 
$(x_i^{[j]},y_i^{[j]})$, $i=0,\ldots,M^{[j]}$, $M^{[j]}\in\mathbb{N}$, i.e. 
\begin{displaymath}
  \gamma^{[j]}(s_i^{[j]}) = \left(\begin{array}{c} x_i^{[j]} \\ y_i^{[j]} \end{array}\right), 
  \qquad i=0,\ldots,M^{[j]}, j=1,\ldots,N.
\end{displaymath} 
Herein, the curves are parametrized with respect to their arc lengths with 
\[s^{[j]}_0:=0,\quad s^{[j]}_{i+1}:=s^{[j]}_{i} + \sqrt{\left(x^{[j]}_{i+1}-x^{[j]}_{i} \right)^2 + \left( y^{[j]}_{i+1}-y^{[j]}_{i}\right)^2 },\quad i=0,\ldots,M^{[j]}-1,\]
and $L^{[j]}:=s_{M^{[j]}}$. 
The initial value problem
\begin{equation}\label{EQ:MPC:Spline}
  (s^{[j]})'(t)=v^{[j]}(t),\qquad s^{[j]}(0)= 0
\end{equation}
describes the motion of the j-th vehicle alongside the spline curve $\gamma^{[j]}$, where $t$ denotes the time and $v^{[j]}(t)$ the velocity of the j-th vehicle at time $t$. 
We assume that we are able to control the velocity $v^{[j]}(t)\in [v^{[j]}_{min},v^{[j]}_{max}]$ of the vehicle, where $v^{[j]}_{min}$ and $v^{[j]}_{max}$ are the minimum and maximum velocity, respectively. 
The position of the j-th vehicle at time $t$ is then given by $\gamma^{[j]}(s^{[j]}(t)) = (x^{[j]}(s^{[j]}(t)),y^{[j]}(s^{[j]}(t)))^\top$. 

In Step 2a of the distributed hierarchical model predictive control algorithm, the j-th vehicle has to solve the 
following optimal control problem on the time horizon $[t_i,t_i + T]$, where $s^{[j]}_\star$ denotes the terminal arc-length 
for vehicle $j$ and $r>0$ is a given security distance. For simplicity we use a ball-shaped constraint for collision avoidance.  

\medskip
{\em
Minimize
\begin{align*}
   \frac{1}{2} ( s^{[j]}(t_i + T) - s^{[j]}_\star )^2
\end{align*}
\indent subject to the constraints 
\begin{align*}
  (s^{[j]})'(t) & =  v^{[j]}(t), &&  t\in [t_i,t_i+T],\\
  v^{[j]}(t) & \in  [v^{[j]}_{min},v^{[j]}_{max}], &&  t\in [t_i,t_i+T],\\
  \| \gamma^{[j]}(s^{[j]}(t)) - \gamma^{[k]}(s^{[k]}(t)) \|^2 & \geq r^2 , & & \forall k\in I_{PR}(j)(t_i), t\in [t_i,t_i+T]. 
\end{align*}}

Instead of controlling the velocity directly it is also possible to control the acceleration of the vehicles, compare (\ref{EQ:2a}), 
or to introduce a delay as in (\ref{EQ:2aa}).

\subsection{Numerical Results}
We present numerical results for the distributed hierarchical model predictive controller in Sections~\ref{sec:Priority}, \ref{sec:OCP-MPC} using 
the full maneuverability of the cars. For a numerical solution it is important to choose suitable parameters. 
Especially the selection of the prediction horizon length $T$ and the control horizon length $\tau$ as well as the weights $\alpha_1,\alpha_2,\alpha_3$ is essential. 
If the prediction horizon is too short the reaction time might be to brief, if it is too large the computational effort is too high. The computational effort can also 
be reduced by the choice of the length of the control horizon, but its size also influences the approximation error. The choice of the weights for the controls 
are linked to the range of the controls and the preferred driving style.

We tested our approach for several everyday scenarios for $N=2$ or $N=3$ cars. All cars are subject to the same 
car model so they have the same dynamics and the same box-constraints as well as the same limits for velocity and 
steering angle. We also assumed that the cars are driving in an equal way, i.e. the objective function 
\eqref{EQ:MPC:OBJ} of each car has the same weights. Furthermore the cars are allowed to occupy the entire road.
Since all cars have the same parameters, we suppress the index $j$ throughout and used the following values: 
\[
\begin{array}{lcllcllcl}
T 		 &=& 2 [s]  & \qquad \tau 	  &=& 0.1 [s] & \qquad          && \\ 
\alpha_1 &=& 1   & \qquad \alpha_2 &=& 1   & \qquad \alpha_3 &=& 10 \\ 
v_{\min} &=& 1 [m/s]  & \qquad v_{\max} &=& 10 [m/s] & \qquad        &&\\ 
a_{\min} &=& -10 [m/s^2] & \qquad a_{\max} &=& 1.5 [m/s^2] & \qquad w_{\max} &=& 0.5 [rad] \\
\ell & = & 4 [m] & \qquad r_x & = & 3.5 [m] & \qquad r_y & = & 2.5 [m]
\end{array} 
\]

\begin{example}[Scenario 1: Avoiding a parking car]
In this scenario we consider two consecutive cars which drive in the same lane and direction as shown in Figure \ref{Scenario_1}. The car in the front is going straight for some time and then it parks in the right lane. The second car is driving behind the first car with the same velocity until the car in the front stops. The parked car is then considered to be an obstacle by the second car and an evasive maneuver is executed. After the second car passed by the first car it changes from the left to right lane again. 

\begin{figure}[h]
  \begin{center}
    \includegraphics[height=7.5cm]{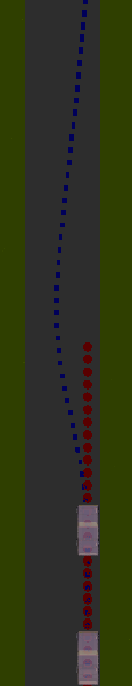}
    \includegraphics[height=7.5cm]{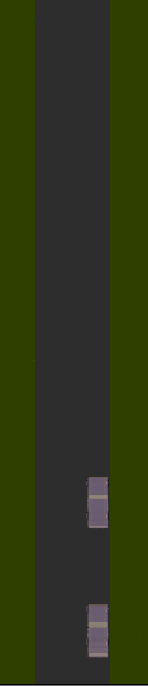}
    \includegraphics[height=7.5cm]{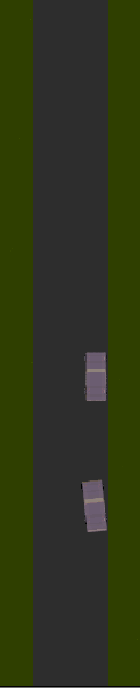}
    \includegraphics[height=7.5cm]{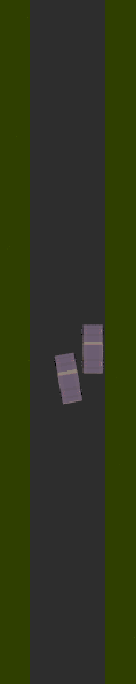}
    \includegraphics[height=7.5cm]{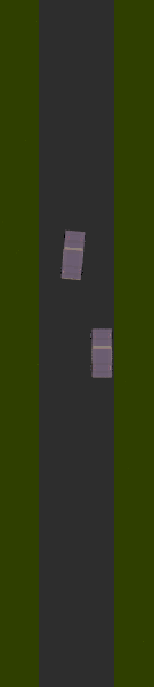}
    \includegraphics[height=7.5cm]{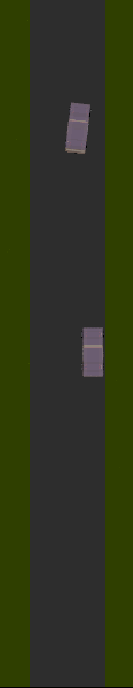}
    \includegraphics[height=7.5cm]{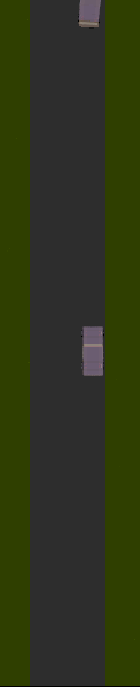}
  \end{center}
  \caption{Scenario 1: Trajectory (left) and snapshots of the motion (pictures on the right).}
  \label{Scenario_1}
\end{figure}

\end{example}

\begin{example}[Scenario 2: Two cars driving through a narrow space.]
This time we examine two cars cars which move towards each other from opposite directions and both cars have to pass through a narrow space in the middle of the road as shown in Figure \ref{Scenario_2}. The car coming from below is closer to the obstacle. Therefore it drives through the narrow space first while the second car slows down and waits until the first car passed through the obstacle. Then the second car accelerates and also moves through the narrow space.

\begin{figure}[h]
  \begin{center}
    \includegraphics[height=7.5cm]{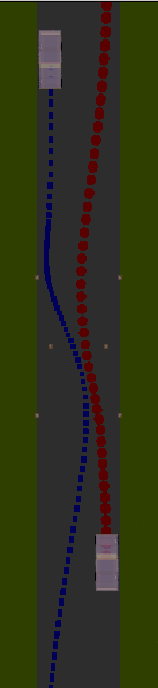}
    \includegraphics[height=7.5cm]{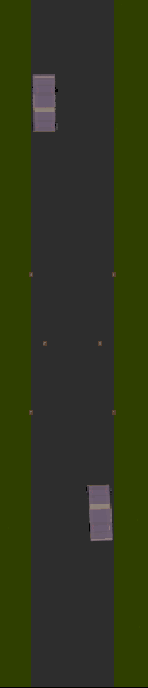}
    \includegraphics[height=7.5cm]{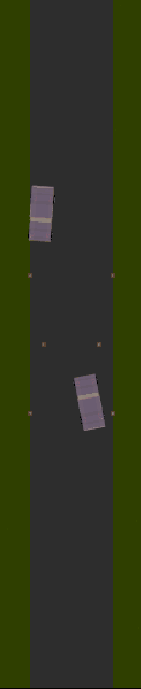}
    \includegraphics[height=7.5cm]{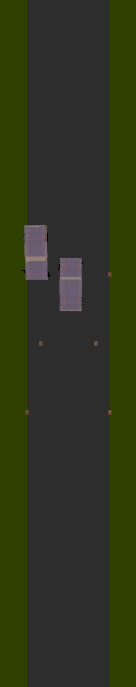}
    \includegraphics[height=7.5cm]{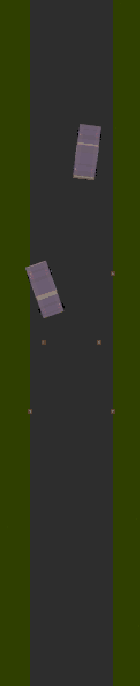}
    \includegraphics[height=7.5cm]{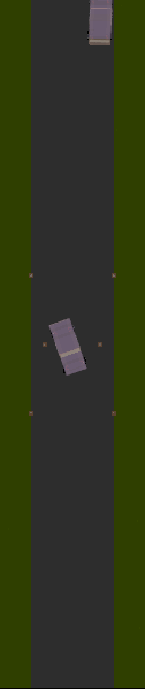}
    \includegraphics[height=7.5cm]{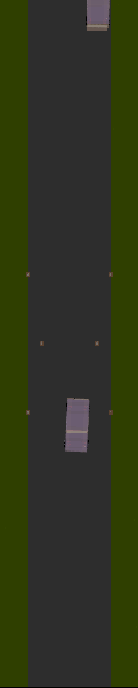}
  \end{center}
  \caption{Scenario 2: Trajectory (left) and snapshots of the motion (pictures on the right).}
  \label{Scenario_2}
\end{figure}
\end{example}

\begin{example}[Scenario 3: Three cars at a intersection]
For the last scenario we consider three cars which cross each other at a intersection, compare Figure \ref{Scenario_3}. In this case, traffic rules apply and the cars on the right of each car have higher priority, i.e. the car coming from below has a higher hierarchy level than the car coming from the left, which has a higher priority than the car coming from above, which has to drive onto the left lane of the road so that the  car in the middle is able to turn left without collision. The car from below is allowed to ignore the other cars and is able to turn left without collision as shown in figure \ref{Screenshots_3}.

\begin{figure}[h]
  \begin{center}  
    \includegraphics[width=0.9\textwidth]{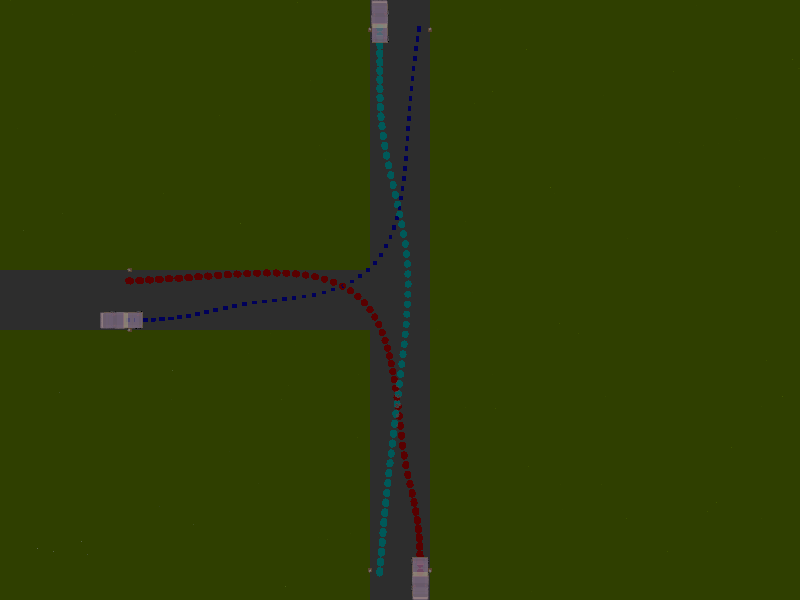}
  \end{center}
  \caption{Scenario 3: Trajectories.}
  \label{Scenario_3}
\end{figure}

\begin{figure}[h]
  \begin{center}
    \includegraphics[width=0.3\textwidth]{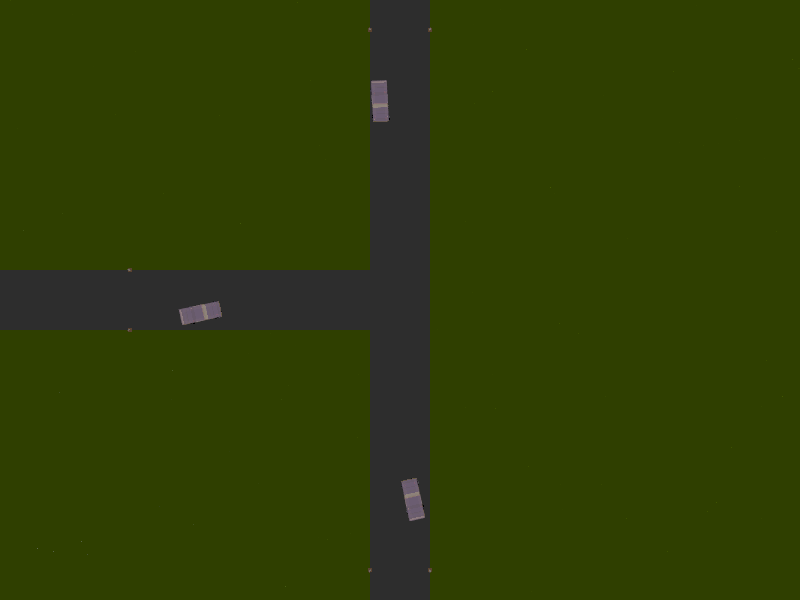}
    \includegraphics[width=0.3\textwidth]{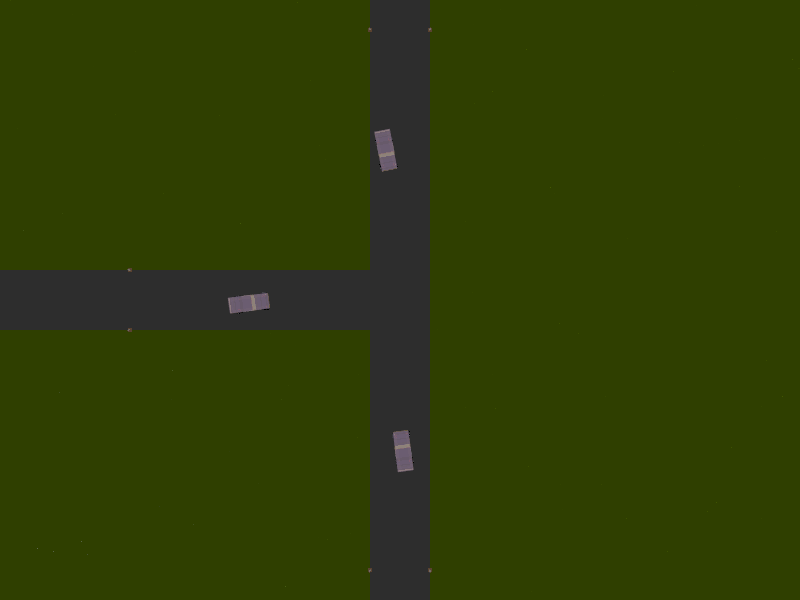}
    \includegraphics[width=0.3\textwidth]{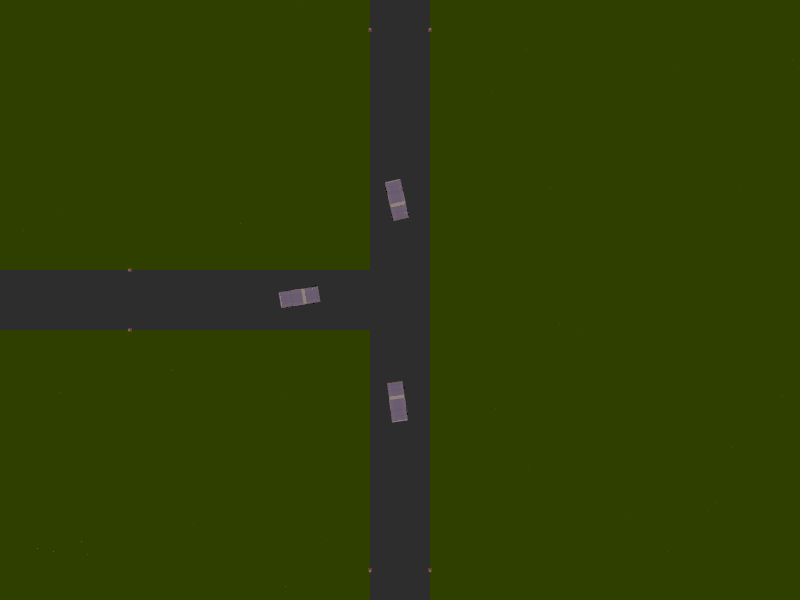}
    \includegraphics[width=0.3\textwidth]{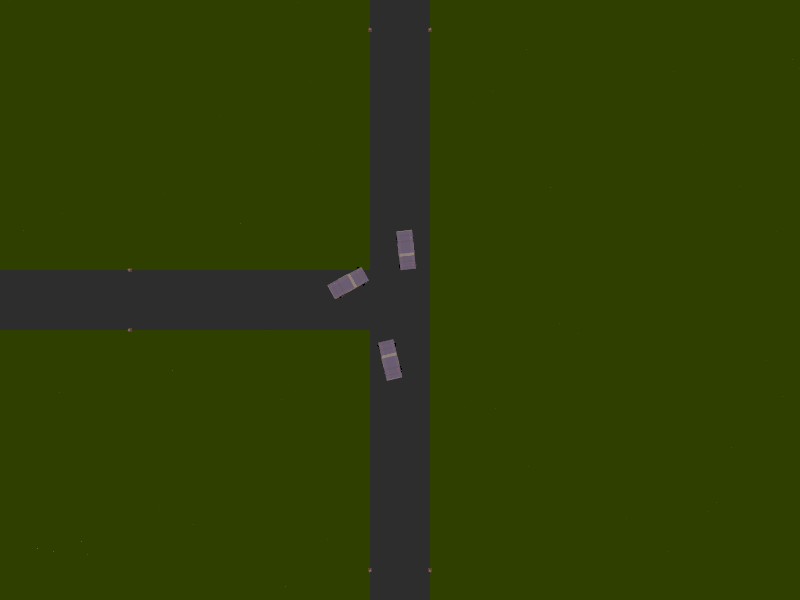}
    \includegraphics[width=0.3\textwidth]{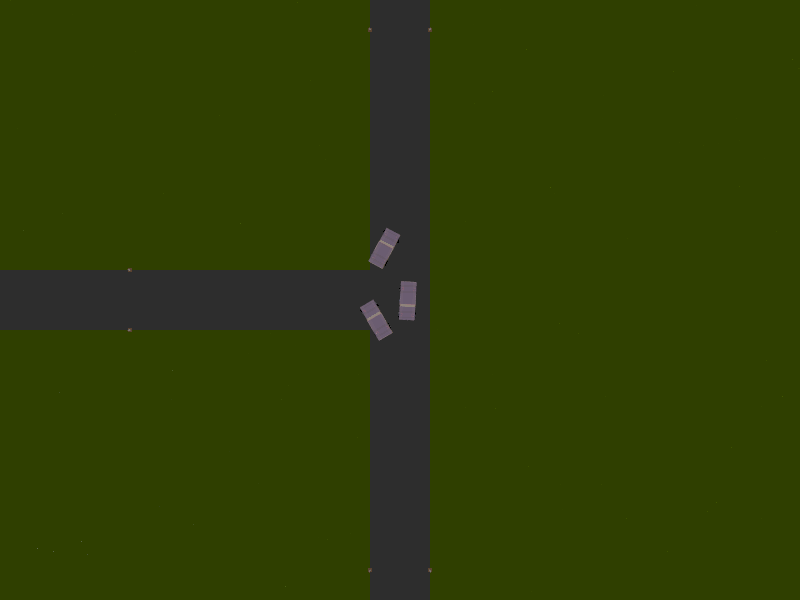}
    \includegraphics[width=0.3\textwidth]{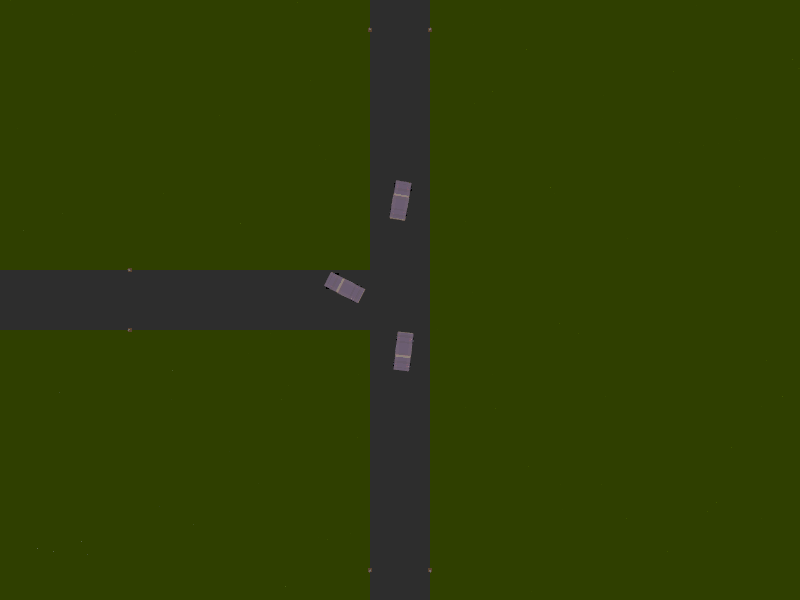}
  \end{center}
  \caption{Snapshots of scenario 3.}
  \label{Screenshots_3}
\end{figure}
\end{example}

\section{A Tracking Controller for Spline Curves}
\label{sec:Tracking}

Once a reference track has been obtained, e.g., by the techniques in Section~\ref{sec:trajectory},
the task is to follow the track with a real car. To this end let the reference track (=desired track) 
be given by a cubic spline curve
\begin{equation}\label{EQ:Spline}
  \gamma(t) := \left(\begin{array}{c} x_d(t) \\ y_d(t) \end{array}\right), \qquad t\in [t_0,t_f], 
\end{equation}
which interpolates the solution obtained by one of the techniques in Section~\ref{sec:trajectory} at 
given grid points within $[t_0,t_f]$. The goal is to design a nonlinear feedback controller 
according to the flatness concept in \cite{Fliess95,Rotella2002,Martin2003}. 
The basic idea is to use inverse kinematics in order to find a feedback law for the control 
inputs of the car. To this end we consider the system of differential equations given by 
(\ref{EQ:1a})-(\ref{EQ:1c}), where we control the velocity $v$ and the steering angle 
$\delta$. We assume that we can measure the outputs 
\begin{displaymath}
  y_1 := x\qquad \mbox{and}\qquad y_2 := y,
\end{displaymath}
i.e. the (x,y)-position of the vehicle. By differentiation we obtain
\begin{eqnarray}
  y_1'  & = & x'  = v \cos \psi, \label{EQ:Out:1a}\\ 
  y_1'' & = & x'' = v' \cos \psi - v \psi'\sin \psi = v' \cos \psi - \frac{v^2}{\ell}\sin \psi \tan \delta,\label{EQ:Out:2a}\\
  y_2'  & = & y'  = v \sin \psi, \label{EQ:Out:1b}\\
  y_2'' & = & y'' = v' \sin \psi + v \psi'\cos \psi = v' \sin \psi + \frac{v^2}{\ell}\cos \psi \tan \delta.\label{EQ:Out:2b}
\end{eqnarray}
Multiplication of \eqref{EQ:Out:1a} by $\cos\psi$ and of \eqref{EQ:Out:1b}
by $\sin\psi$, adding both equations, exploiting $\psi = \arctan(y'/x')$ and solving for the control $v$ yields 
\begin{eqnarray*}
  v = V(y_1,y_2) & := & y_1'\cos\left( \arctan\left(\frac{y_2'}{y_1'} \right)\right)  + y_2'\sin\left( \arctan\left(\frac{y_2'}{y_1'} \right)\right) \\
                & = & \sqrt{ (y_1')^2 + (y_2')^2}.
\end{eqnarray*}
Likewise, multiplication of \eqref{EQ:Out:2a} by $\sin\psi$ and of \eqref{EQ:Out:2b} 
by $\cos\psi$, subtracting equations, exploiting $\psi = \arctan(y'/x')$, $v = \sqrt{(y_1')^2 +(y_2')^2}$, and solving for the control $\delta$ yields 
\begin{eqnarray*}
  \delta = \Delta(y_1,y_2) & := & \arctan\left(\frac{\ell \left( y_2''\cos\left( \arctan\left(\frac{y_2'}{y_1'} \right)\right)  - y_1''\sin\left( \arctan\left(\frac{y_2'}{y_1'}\right)\right)\right)   }{v^2} \right) \\
  & = & \arctan\left(\frac{\ell \left( y_2'' y_1' - y_1'' y_2'\right)}{\left(\left(y_{1}'\right)^2 + \left(y_{2}'\right)^2\right)^\frac{3}{2}} \right).
\end{eqnarray*}
If we introduce the reference coordinates $x_d$ and $y_d$ (=desired outputs) from (\ref{EQ:Spline}) into these formula, then the corresponding controls 
$v_d = V(x_d,y_d)$ and $\delta_d = \Delta(x_d,y_d)$ would 
track the reference input provided the initial value is consistent with the reference trajectory. However, in practice there will be deviations due to 
modeling errors or disturbances. Hence, we need a feedback control law that is capable of taking deviations from the reference input into account. Such a feedback control law reads as follows, compare \cite{Rotella2002} for general flat systems: 
\begin{eqnarray*}
  & & K_v(y_1,y_2,y_{1,d},y_{2,d}) := \sqrt{\left( y_{1,d}' - k_1\left(y_{1} - y_{1,d}\right) \right)^2 + \left( y_{2,d}' - k_2\left( y_2 - y_{2,d}\right) \right)^2}, 
\end{eqnarray*}
and 
\begin{eqnarray*}
  & & K_\delta(y_1,y_2,y_{1,d},y_{2,d}) := \\
  & & \qquad \arctan\left(\ell\,\frac{\left( y_{2,d}'' - k_5\left(y_{2}' - y_{2,d}'\right) - k_6\left( y_2 - y_{2,d}\right) \right)\,y_{1}'}{\left(\left(y_{1}'\right)^2 + \left(y_{2}'\right)^2\right)^\frac{3}{2}} \right. \\
  & &\qquad\qquad\qquad \left. -\ell  \frac{\left( y_{1,d}'' - k_3\left(y_{1}' - y_{1,d}'\right) - k_4\left( y_1 - y_{1,d}\right) \right)\,y_{2}'}{\left(\left(y_{1}'\right)^2 + \left(y_{2}'\right)^2\right)^\frac{3}{2}} \right) .
\end{eqnarray*}
Herein, $k_1,\ldots,k_6$ are constants that influence the response time. 
$y_1 \ (=x)$ and $y_2\ (=y)$ are the actual measurements of the car's midpoint position of the rear axle and $y_{1,d}\ (=x_d)$ and $y_{2,d} \ (=y_d)$ are the reference coordinates, respectively. In addition, the derivatives $y_1'=x'$ and $y_2'=y'$ need to be estimated as well, for instance by finite difference approximations using the position measurements. 

If we insert this feedback control law into our system, we obtain the closed loop system
\begin{eqnarray}
  x'(t) & =  K_v(x(t),y(t),x_d(t),y_d(t))\cos\psi(t), \label{EQ:Feedback:1a}\\
  y'(t) & =  K_v(x(t),y(t),x_d(t),y_d(t))\sin\psi(t), \label{EQ:Feedback:1b}\\
  \psi'(t) & =\frac{K_v(x(t),y(t),x_d(t),y_d(t))}{\ell} \tan K_\delta(x(t),y(t),x_d(t),y_d(t)).\label{EQ:Feedback:1c}
\end{eqnarray}
To study the stability behavior we assume $k_1=k_2,k_3=k_5,k_4=k_6$ and linearize the right side of the system with respect to the reference trajectory, which gives us the time variant matrix
\begin{equation*}
A =
\begin{pmatrix}
-\frac{k_1\,\left(x_d'\right)^2}{v_d} & -\frac{k_1\,x_d'\,y_d'}{v_d} & -y_d'\\ 
-\frac{k_1\,x_d'\,y_d'}{v_d}  & -\frac{k_1\,\left(y_d'\right)^2}{v_d} & x_d'\\ 
\frac{k_1\,\left(y_d''\,\left(x_d'\right)^2 - x_d''\,x_d'\,y_d'\right) - k_3\,y_d'\,v_d^2}{v_d^4} & \frac{k_1\left(y_d''\,x_d'\,y_d' - x_d''\,\left(y_d'\right)^2\right) - k_3\,y_d'\,v_d^2}{v_d^4} & -k_2
\end{pmatrix} ,
\end{equation*}
with $v_d:=\sqrt{\left(x_d'\right)^2 + \left(y_d'\right)^2}$, where we suppressed the argument $t$ for notational simplicity. The characteristic polynomial of the matrix is
\begin{equation*}
\det(\lambda\,I_4 - A)=\left(\lambda^2+k_3\lambda+k_4\right)\left(\lambda+k_1\right).
\end{equation*}
It follows, if $\lambda+k_1$ and $\lambda^2+k_3\lambda+k_4$ are Hurwitz the linearized system is asymptotic stable and thereby the closed loop nonlinear system is locally asymptotic stable.

Measurement errors can be modeled by introducing white noise $(x_{wn},y_{wn})^\top$ and $(x'_{wn},y'_{wn})^\top$ into the feedback control law, i.e. by changing $x$ and $y$ to $x-x_{wn}$ and 
$y-y_{wn}$ and by changing $x'$ and $y'$ to $x'-x'_{wn}$ and $y'-y'_{wn}$.

\begin{example}[Tracking controller]
  To test the constructed feedback controller we consider to have a constant reference velocity $v_d=11.5\,[\frac{m}{s}]$ and $\ell=2.8\,[m]$ to be the distance from rear axle to front axle. For the parameters in the closed loop system we chose $k_1=k_2=1$, $k_3=k_4=k_5=k_6=2$. 
  For the system with white noise we chose a random number generator with values in the interval $[-10,10)$ (in meters) for the $(x,y)$-position and in the interval 
    $[-2,2]$ (in meters per second) for the velocities $(x',y')$. 
    In practice it is possible to measure the actual position with a tolerance of below one millimeter, but for illustration purposes we use 
    a higher tolerance to test the tracking ability of the controller. For the system with an offset we changed the initial 
    state from $(-26,-1)^\top$ to $(-16,9)^\top$. In all cases the control values of the feedback control law are projected 
    into the feasible control set $[0,50]$ for the velocity $v$ and to $[-\pi/6,\pi/6]$ for the steering angle $\delta$. 

    Figure~\ref{Fig:EX1a} 
    shows the simulations with the feedback 
    controlled system (\ref{EQ:Feedback:1a}) - (\ref{EQ:Feedback:1c}) for a given track at a sampling rate of 20 Hz. For better visibility only every 5th data point 
    is plotted in Figure~\ref{Fig:EX1a}. 
    In all cases the controller is able to accurately follow the track. 
    
    \begin{figure}[h]
      \begin{center}
        \includegraphics[height=7cm,angle=-90]{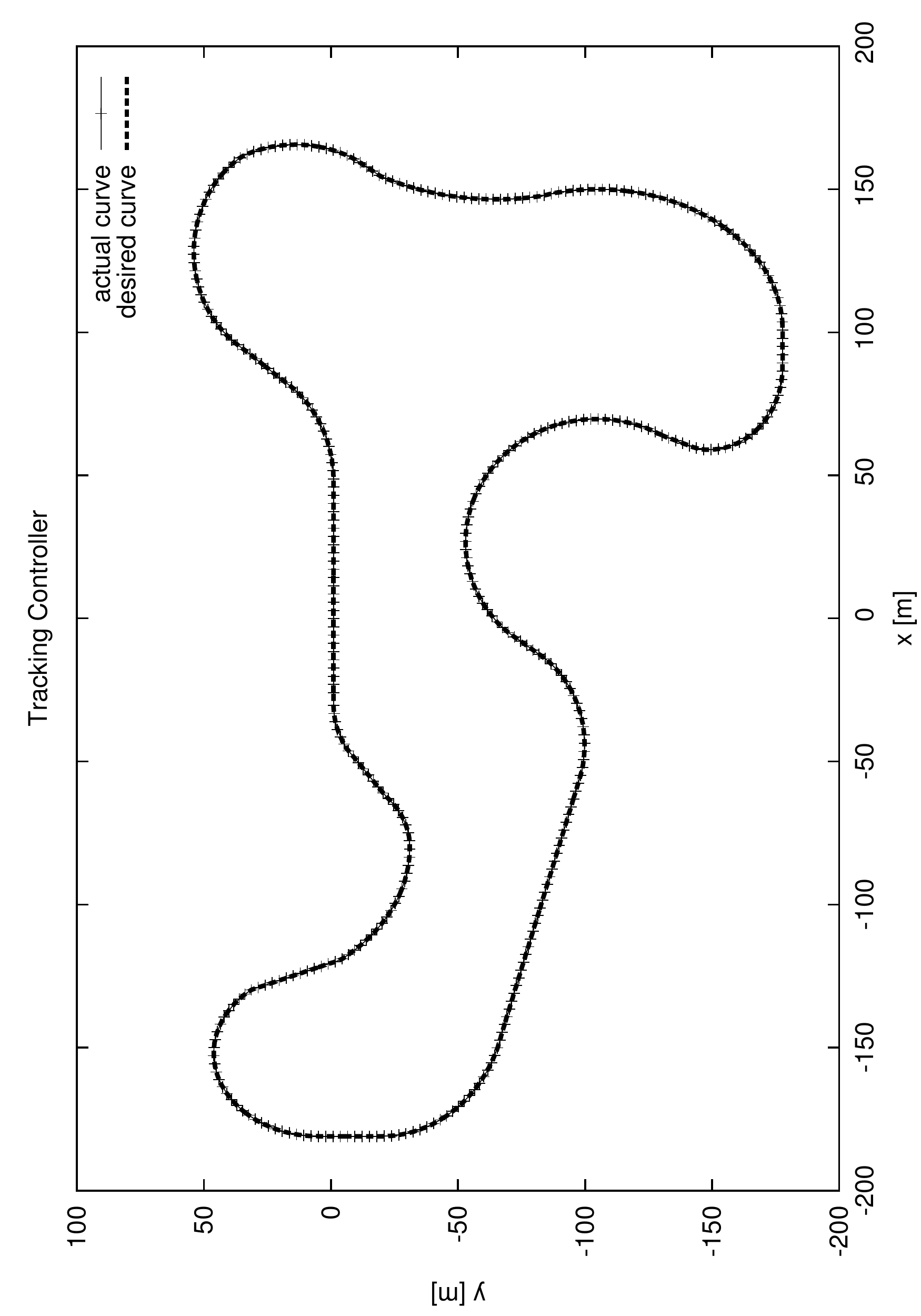}
        \includegraphics[height=7cm,angle=-90]{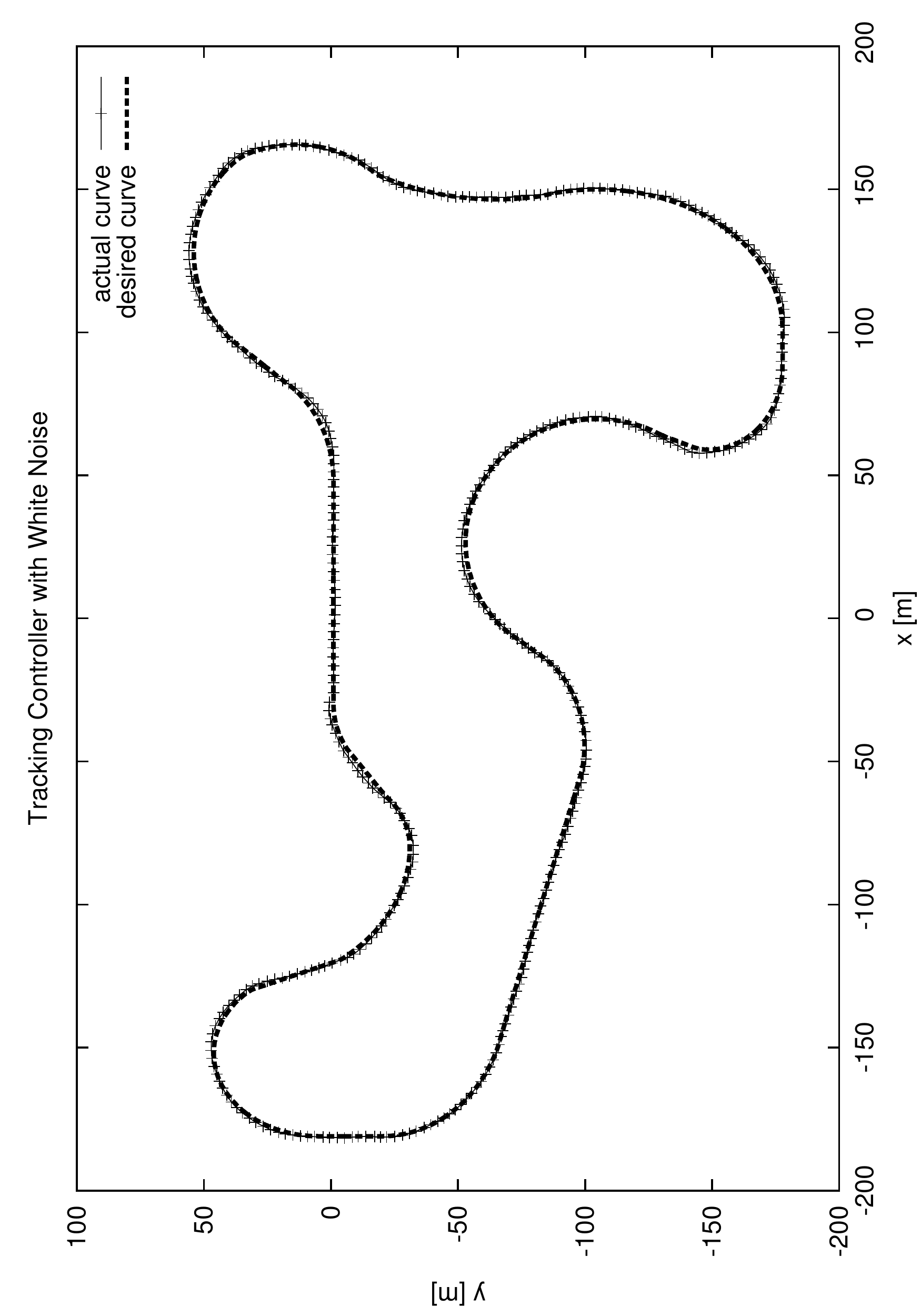}
        \includegraphics[height=7cm,angle=-90]{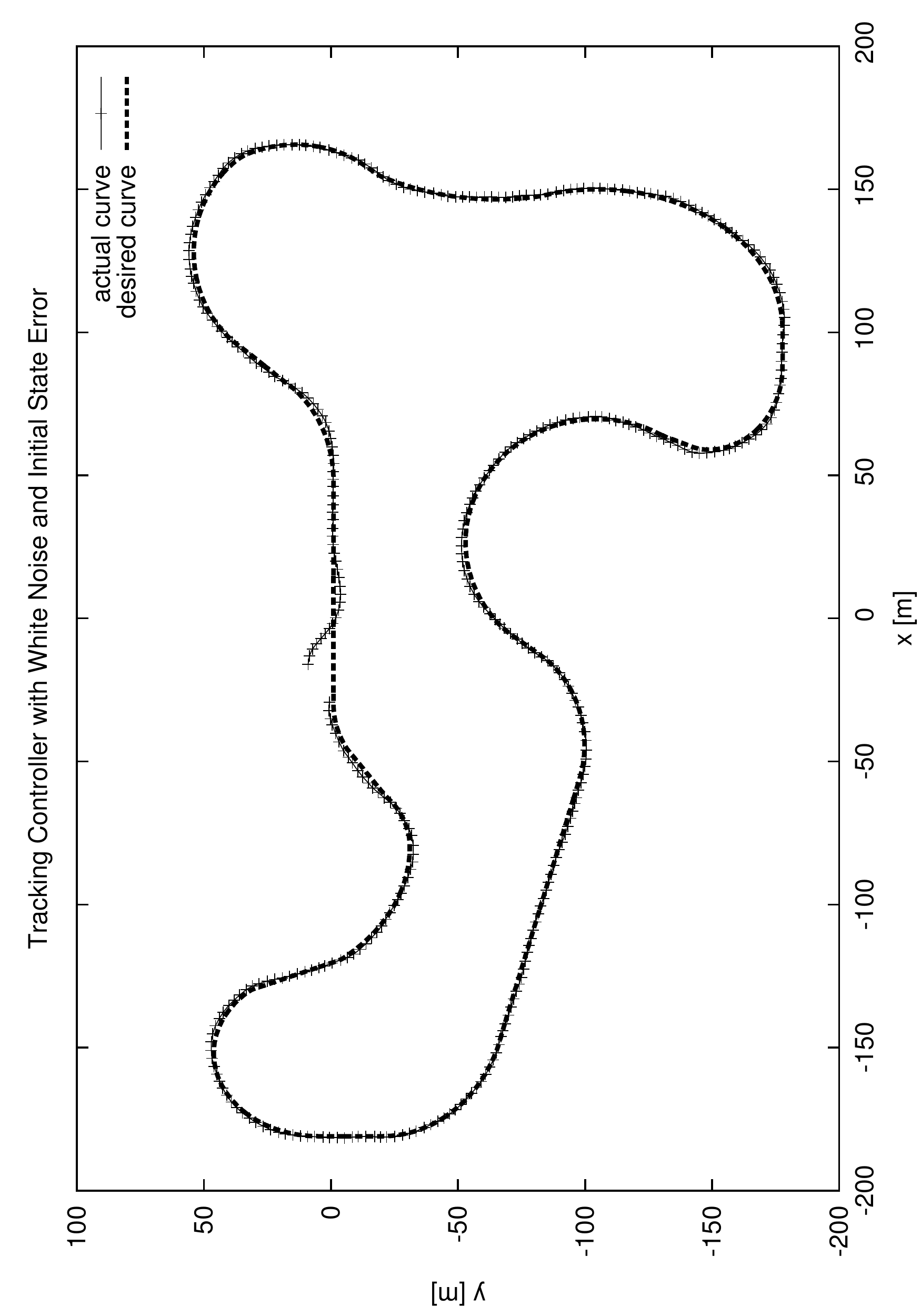}
      \end{center}
      \captionof{figure}{Tracking controller output: Trajectory of car with initial position on the track (top), with white noise (middle), and with white noise and initial state offset (bottom).}
      \label{Fig:EX1a}
    \end{figure}
\end{example}

\begin{acknowledgement}
  This material is based upon work supported by the Air
   Force Office of Scientific Research, Air Force Materiel Command, USAF, under Award
   No, FA9550-14-1-0298. 

\end{acknowledgement}
%
%
%
%
\bibliographystyle{spmpsci}

\end{document}